\documentclass{article}

\title{Ricci Flow method in the existence problem of the K$\ddot{a}$hler-Einstein metrics}

\author{LIU CHAO}
\date{}
\begin{document}
	\maketitle
	\section{Introduction}

This note  illustrates the Ricci flow method based on the Cao.H.D's paper[1] and Yau.S.T's paper[4], and tries to explain the method in detail, especially in some calculations.  Jian Song and Weinkove's note[9] used some other estimates to obtain the result, this paper will explain some of their estimates as well. This note was a seminar lecture note in 2022 summer when the author was giving lectures on the geometry analysis seminar reasearching the  Ricci flow method. The part of the imporatnt zero order estimate is going to be added in a few days. The level of the author is limited, if there are any errors, please do not hesitate to advise.  Any comments will be grateful.
    \section{The Ricci Flow equation}
    
    Let M be a compact K$\ddot{a}$hler manifold with the K$\ddot{a}$hler dimension n  and the K$\ddot{a}$hler metric 
    $\mathit{d}$s = $\mathit{g_{i\bar{j}}}\mathit{d}\mathit{z^{i}}\wedge\mathit{d}\mathit{\bar{z}^{j}}$. We will use the Einstein summation through the whole article.

    Let
  $\mathit{R_{i\bar{j}}}$ = -$\frac{\partial^{2}}{\partial\mathit{z^{i}}\partial\mathit{\bar{z}^{j}}}$log det($\mathit{g_{i\bar{j}}}$) 
     be the Ricci curvatureof the corresponding metric, and the Ricci Form, i.e. the (1,1) tensor $\frac{\sqrt{-1}}{2\pi}\mathit{R_{i\bar{j}}}\mathit{d}\mathit{z^{i}}\wedge\mathit{d}\mathit{\bar{z}^{j}}$ 
     is closed so we can define the cohomology class of it which is the first chern class C$_{1}$(M) of M.

    We consider the Ricci Flow equation
     $\frac{\partial{g_{ij}}}{\partial t}$ = -2$\mathit{R_{i\bar{j}}}$  + $\frac{2}{3}$r$\mathit{g_{ij}}$ 
     and its complex version:
    
    $\frac{\partial\tilde{\mathit{g_{i\bar{j}}}}}{\partial t}$ = -$\tilde{\mathit{R_{i\bar{j}}}}$ + $\mathit{T_{i\bar{j}}}$ ,  $\tilde{\mathit{g_{i\bar{j}}}}$ = $\mathit{g_{i\bar{j}}}$ at t = 0         (2.1)

    where $\tilde{\mathit{R_{i\bar{j}}}}$
     is the Ricci tensor of the
      $\tilde{\mathit{g_{i\bar{j}}}}$,  $\mathit{T_{i\bar{j}}}$ is a representation of the first Chern class C$_{1}$(M), actually we can choose any representation satisfying our requirement, such as a K$\ddot{a}$hler-Einstein metric.
    
    Then our idea is:first prove the solution exists all the time, thenprove when t goes infinitly, 
     $\tilde{\mathit{g_{i\bar{j}}}}$
      converges to a definite 
       $\tilde{\mathit{g_{i\bar{j}}}}$($\infty$)
        and hence 
        $\frac{\partial\tilde{\mathit{g_{i\bar{j}}}}}{\partial t}$
         converges to 0, then we get 
         -$\tilde{\mathit{R_{i\bar{j}}}}$ + $\mathit{T_{i\bar{j}}}$ = 0, 
         then we get $\mathit{T_{i\bar{j}}}$ 
         will be the Ricci tensor of
          $\tilde{\mathit{g_{i\bar{j}}}}$($\infty$)
           so the $\tilde{\mathit{g_{i\bar{j}}}}$($\infty$) is the metric we want.
    
    However,  the equation (2.1) is too abstract to solve, therefore, we replace it.Due to
     $\frac{\sqrt{-1}}{2\pi}\mathit{T_{i\bar{j}}}\mathit{d}\mathit{z^{i}}\wedge\mathit{d}\mathit{\bar{z}^{j}}$ 
      is in the class of first chern class of M, so is

     $\frac{\sqrt{-1}}{2\pi}\mathit{R_{i\bar{j}}}\mathit{d}\mathit{z^{i}}\wedge\mathit{d}\mathit{\bar{z}^{j}}$ ,  where by the Hodge theory and the $\partial \bar{\partial}$ Lemma, the cohomology group $H^{1,1}_{\bar{\partial}}$(M,R) = $\frac{\{\bar{\partial} \;closed \;real \;(1,1) forms\;\}}{Im\partial \bar{\partial}}$,
     so they deviate by a term which in the kernel of the boundary operator $\mathit{d} $, i.e.
     
      $\frac{\partial^{2}f}{\partial\mathit{z^{i}}\partial\mathit{\bar{z}^{j}}}$:
       $\mathit{T_{i\bar{j}}}$ - $\mathit{R_{i\bar{j}}}$ =  $\frac{\partial^{2}f}{\partial\mathit{z^{i}}\partial\mathit{\bar{z}^{j}}}$. So we take t = 0 in the equation, similary, we can assume
       
        $\tilde{\mathit{g_{i\bar{j}}}}$ - $\mathit{g_{i\bar{j}}}$ = $\frac{\partial^{2}u}{\partial\mathit{z^{i}}\partial\mathit{\bar{z}^{j}}}$, for u$\in$C$^{\infty}(M\times[0,T))$ , 
        T is nonegative and less or equal to infinity,in order to satisfy the initial condition, u(0) = 0, 
        we caculate the new scaler equation replacing  $\tilde{\mathit{g_{i\bar{j}}}}$ by 
        $\tilde{\mathit{g_{i\bar{j}}}}$ = $\mathit{g_{i\bar{j}}}$ + $\frac{\partial^{2}u}{\partial\mathit{z^{i}}\partial\mathit{\bar{z}^{j}}}$:

        The left hand side is:
        $\frac{\mathit{g_{i\bar{j}}} + \frac{\partial^{2}u}{\partial\mathit{z^{i}}\partial\mathit{\bar{z}^{j}}}}{\partial t}$=
        0 + $\partial\frac{ \frac{\partial^{3}u}{\partial\mathit{z^{i}}\partial\mathit{\bar{z}^{j}}}}{\partial t}$ = 
        $\frac{\partial^{2}}{\partial\mathit{z^{i}}\partial\mathit{\bar{z}^{j}}}$($\frac{\partial u}{\partial t}$)

         The right hand side is : 
          -$\tilde{\mathit{R_{i\bar{j}}}}$ + $\mathit{R_{i\bar{j}}}$ + $\frac{\partial^{2}f}{\partial\mathit{z^{i}}\partial\mathit{\bar{z}^{j}}}$

           while we can calculate the Ricci tensor explicitly in  terms of $\mathit{g_{i\bar{j}}}$ and $\tilde{\mathit{g_{i\bar{j}}}}$:
           
            $\mathit{R_{i\bar{j}}}$ = -$\frac{\partial^{2}}{\partial\mathit{z^{i}}\partial\mathit{\bar{z}^{j}}}$ log det($\mathit{g_{i\bar{j}}}$)

             so the right hand side is euqal to :

             -(-$\frac{\partial^{2}}{\partial\mathit{z^{i}}\partial\mathit{\bar{z}^{j}}}$ log det($\tilde{\mathit{g_{i\bar{j}}}}$))- $\frac{\partial^{2}}{\partial\mathit{z^{i}}\partial\mathit{\bar{z}^{j}}}$ log det($\mathit{g_{i\bar{j}}}$) + $\frac{\partial^{2}f}{\partial\mathit{z^{i}}\partial\mathit{\bar{z}^{j}}}$

             =$\frac{\partial^{2}}{\partial\mathit{z^{i}}\partial\mathit{\bar{z}^{j}}}$((log det($\mathit{g_{i\bar{j}}}$ + $\frac{\partial^{2}u}{\partial\mathit{z^{i}}\partial\mathit{\bar{z}^{j}}}$ ) - log det($\mathit{g_{i\bar{j}}}$ )) + $\frac{\partial^{2}f}{\partial\mathit{z^{i}}\partial\mathit{\bar{z}^{j}}}$
             
             So we get:

                $\frac{\partial^{2}}{\partial\mathit{z^{i}}\partial\mathit{\bar{z}^{j}}}$($\frac{\partial u}{\partial t}$) =   
                $\frac{\partial^{2}}{\partial\mathit{z^{i}}\partial\mathit{\bar{z}^{j}}}$((log det($\mathit{g_{i\bar{j}}}$ + $\frac{\partial^{2}u}{\partial\mathit{z^{i}}\partial\mathit{\bar{z}^{j}}}$ ) - log det($\mathit{g_{i\bar{j}}}$ )) + $\frac{\partial^{2}f}{\partial\mathit{z^{i}}\partial\mathit{\bar{z}^{j}}}$

                Finally, 
                 $\frac{\partial^{2}}{\partial\mathit{z^{i}}\partial\mathit{\bar{z}^{j}}}$($\frac{\partial u}{\partial t}$ -
                 log det($\mathit{g_{i\bar{j}}}$ + $\frac{\partial^{2}u}{\partial\mathit{z^{i}}\partial\mathit{\bar{z}^{j}}}$ ) + log det($\mathit{g_{i\bar{j}}}$ ) - f) = 0,  
                  the function
                  
                  ($\frac{\partial u}{\partial t}$ -
                  log det($\mathit{g_{i\bar{j}}}$ + $\frac{\partial^{2}u}{\partial\mathit{z^{i}}\partial\mathit{\bar{z}^{j}}}$ ) + log det($\mathit{g_{i\bar{j}}}$ ) - f) satisfies this uniform elliptic equation $\Delta F=0$, then by the strong maximum principle, the maximum and the minimum attain on the boundry, but M is compact,so only non-boundry points exist, then this function can only be a constant relative to $\partial\bar{\partial}$,
                  this means the term in the partials is equal to a smooth function relative to t:

                 $\frac{\partial u}{\partial t}$ =
                 log det($\mathit{g_{i\bar{j}}}$ + $\frac{\partial^{2}u}{\partial\mathit{z^{i}}\partial\mathit{\bar{z}^{j}}}$ ) - log det($\mathit{g_{i\bar{j}}}$ ) + f + $\phi(t)$

                And the $\phi(t)$ should satisfy the compatibility condition:
                
                $\int_{M} e^{\frac{\partial u}{\partial t} - f}\mathit{d}$$\tilde{V}$ = $e^{\varphi(t)}$Vol(M)
                
                Actually, this compatibility condition wants to say that the volume of the compact manifold stays invariant during the deformation of the metric $\tilde{\mathit{g}}^{i\bar{j}}$(t).
                From the equation above, 
                
                $\frac{\partial u}{\partial t}$ - f =
                log det($\mathit{g_{i\bar{j}}}$ + $\frac{\partial^{2}u}{\partial\mathit{z^{i}}\partial\mathit{\bar{z}^{j}}}$ ) - log det($\mathit{g_{i\bar{j}}}$ ) + $\phi(t)$
                
                we take the exponetial of the both sides and integrate them on M:
                
                $\int_{M}e^{\frac{\partial u}{\partial t} -f }\mathit{d}V$ = 
                $\int_{M}e^{log det(\mathit{g_{i\bar{j}}} + \frac{\partial^{2}u}{\partial\mathit{z^{i}}\partial\mathit{\bar{z}^{j}}}) - log det(\mathit{g_{i\bar{j}}} ) + \phi(t)}\mathit{d}V$
                
                 $\int_{M}e^{\frac{\partial u}{\partial t} -f }\mathit{d}V$ = 
                 $e^{\varphi(t)}$ $\int_{M} \frac{logdet(\tilde{\mathit{g_{i\bar{j}}}})}{logdet\mathit{g_{i\bar{j}}}}\mathit{d}V$
                
                While                     $\mathit{d}$V = det($\mathit{g_{i\bar{j}}}$)  $\wedge^{n}_{i = 1}$ ( $\frac{\sqrt{-1}}{2}$ $\mathit{d}z^{i}$ $\wedge$  $\mathit{d}\bar{z}^{j}$) = $\frac{\omega^n}{n!}$

                $\mathit{d}$$\tilde{V}$ = det( $\tilde{\mathit{g_{i\bar{j}}}}$)  $\wedge^{n}_{i = 1}$ ( $\frac{\sqrt{-1}}{2}$ $\mathit{d}z^{i}$ $\wedge$  $\mathit{d}\bar{z}^{j}$)
                 = $\frac{\tilde{\omega} ^n}{n!}$
                
                so the above equation changes to:
                
                $\int_{M}e^{\frac{\partial u}{\partial t} -f }\mathit{d}V$ = 
                $e^{\varphi(t)}\int_{M} \mathit{d}\tilde{V}$  to make sure the volume of M stay unchanged, there is $\int_{M}e^{\frac{\partial u}{\partial t} -f }\mathit{d}V$ = 
                $e^{\varphi(t)}$ $\int_{M} \mathit{d}\tilde{V}$=$e^{\varphi(t)}$ $\int_{M} \mathit{d}V$=$e^{\varphi(t)}$Vol(M).

                To prove the long time existence of the solutoin of this parabolic equation, generally, we need to give up to third order estimate. In the proof of the existence and the uniform convergence, we will use these estimates.

                \section{Existence for the solution in all time}
                
                Actually , the goal equation is given initial value:

                 $\frac{\partial u}{\partial t}$ =
                log det($\mathit{g_{i\bar{j}}}$ + $\frac{\partial^{2}u}{\partial\mathit{z^{i}}\partial\mathit{\bar{z}^{j}}}$ ) - log det($\mathit{g_{i\bar{j}}}$ ) + f \   \    \   \  \  \  \  \  (3.1)

                u(x,t) = 0 when t = 0. And by the  initial assumption, the solution exists in the interval (0,T] and the Kahler metric $\tilde{\mathit{g_{i\bar{j}}}}$ = $\mathit{g_{i\bar{j}}}$ + $\frac{\partial^{2}u}{\partial\mathit{z^{i}}\partial\mathit{\bar{z}^{j}}}$ is positive definite so it's a K$\ddot{a}$hler metric for any t$\in$(0,T].

                To prove the estimation, we need some notations. We differentiate the equation (3.1):

                $\frac{\partial}{\partial t}$($\frac{\partial u}{\partial t}$) =  $\frac{\partial}{\partial t}$(log det($\mathit{g_{i\bar{j}}}$ + $\frac{\partial^{2}u}{\partial\mathit{z^{i}}\partial\mathit{\bar{z}^{j}}}$ )) , for  $\frac{\partial}{\partial t}$(f) = 0 ,and $\frac{\partial}{\partial t}$(log det($\mathit{g_{i\bar{j}}}$) = 0
                then  
                
                $\frac{\partial}{\partial t}$(log det($\mathit{g_{i\bar{j}}}$ + $\frac{\partial^{2}u}{\partial\mathit{z^{i}}\partial\mathit{\bar{z}^{j}}}))$ = 
                $\frac{\frac{\partial}{\partial t}{det( \tilde{\mathit{g_{i\bar{j}}}})}}{det( \tilde{\mathit{g_{i\bar{j}}}})}$
                 
                 = $\frac{1}{det( \tilde{\mathit{g_{i\bar{j}}}})}$$ \tilde{\mathit{g}}^{i\bar{j}}$ $\frac{\partial }{\partial t}\tilde{\mathit{g_{i\bar{j}}}}$$det( \tilde{\mathit{g_{i\bar{j}}}})$
                 
                  = $ \tilde{\mathit{g}}^{i\bar{j}}$ $\frac{\partial }{\partial t}$ $(\frac{\partial^{2}u}{\partial\mathit{z^{i}}\partial\mathit{\bar{z}^{j}}})$ 
                  where ${\tilde\mathit{g}^{i\bar{j}}}$ is the inverse of $\tilde{\mathit{g_{i\bar{j}}}}$.
                  It can not be ignored that the $ \tilde{\mathit{g}}^{i\bar{j}}$ $\frac{\partial }{\partial t}$ $(\frac{\partial^{2}u}{\partial\mathit{z^{i}}\partial\mathit{\bar{z}^{j}}})$ is a Einstein summation, the ${\tilde\mathit{g}^{i\bar{j}}}$ is actually the $((\tilde{\mathit{g_{i\bar{j}}}})^{-1})^{i\bar{j}}$.

               So we get $\frac{\partial}{\partial t}$($\frac{\partial u}{\partial t}$) = ${\tilde\mathit{g}^{i\bar{j}}}$ $\frac{\partial^{2}}{\partial\mathit{z^{i}}\partial\mathit{\bar{z}^{j}}}$($\frac{\partial u}{\partial t}$) , this is a parabolic equation ,so by the maximum principle, at the time t = 0, we can consider this manifold with the initial metric is the bounder of the domain M$\times$[0,T). So we have 
               
               max$_{M}$$\mid$ $\frac{\partial u}{\partial t}$$\mid$ $\le$ max$_{t = 0}$$\mid$$\frac{\partial u}{\partial t}$$\mid$ =  max$_{M}$$\mid$f$\mid$.
                
                Let $\tilde{\Delta} $ = $\tilde{\mathit{g_{i\bar{j}}}}$$\frac{\partial^{2}}{\partial\mathit{z^{i}}\partial\mathit{\bar{z}^{j}}}$ be the normalized Laplace of $\tilde{\mathit{g_{i\bar{j}}}}$, and $\Delta $ =   $\mathit{g_{i\bar{j}}}$$\frac{\partial^{2}}{\partial\mathit{z^{i}}\partial\mathit{\bar{z}^{j}}}$ be the normalized  Laplace of the  $\mathit{g_{i\bar{j}}}$, $\diamond$ = $\tilde{\Delta} $ - $\frac{\partial}{\partial t}$.
                
                
                $\mathbf{Zero\; order\; estimate}$

                Let v = u - $\frac{1}{Vol(M)}$$\int_{M}u\mathit{d}V$ be the normalized of u such that $\int_{M} v = 0$ then  v satisfies the equation:
                
                  det($\mathit{g_{i\bar{j}}}$ + $\frac{\partial^{2} v}{\partial\mathit{z^{i}}\partial\mathit{\bar{z}^{j}}}$) $(det(\mathit{g_{i\bar{j}}}))^{-1}$ = $e^{F}$ 
                  
                  where F = $\frac{\partial u}{\partial t} - f$
                  in the paper of Yau, so we can directly use the calculation in Yau's work, we get the following Lemma:

                $\mathbf{Lemma \;1}$
                
               $ sup_{M\times [0,T)}$ v$\le$ $C_{2}$,   $ sup_{M \times [0,T)}$ $\int_{M}$ $\mid v \mid \mathit{d} V $ $\le$ $C_{3}$  for $C_{2}$ and $C_{3}$ constants
             
          $ \mathbf{Remark}$
             We are here to give a sketch proof of this zero order estimate following Yau's method. We consider the Green's funciton G(p,q) of the normalized Laplace operator of $\mathit{g_{i\bar{j}}}$ ,$\Delta$ on M, and let K be a constant which depends only on M such that G(p,q) + K $\ge$ 0. Then we find that 
             
             $\Delta v$ = $\Delta u$ so
             $\tilde{\mathit{g_{i\bar{j}}}}$ =  $\mathit{g_{i\bar{j}}}$ + $\frac{\partial^{2}v}{\partial\mathit{z^{i}}\partial\mathit{\bar{z}^{j}}}$.
             and $\int_{M} v$ = 0, then by Yau's work it shows
              v(p) = $-\int_{M}(G(p,q) + K)\mathit{d}q$ and $sup_{M}$ $\le$ n $sup_{p\in M}\int_{M}(G(p,q) + K)\mathit{d}q$,
              while the right hand side depends only on M, hence we can consider it as a constant $C_{2}$, and we proceed by this estimate we have 
              
              $\int_{M} \mid v \mid$ $\le$ $\int_{M} \mid sup_{M}v - v \mid $ + $\int_{M} \mid sup_{M}v \mid$ 
              
              $\le$ ($sup_{M} v$ )Vol(M) - $\int_{M} v$ + ($sup_{M} v$ )Vol(M) 
              
              $\le$ 2n Vol(M) $sup_{p\in M}$ $\int_{M} (G(p,q) + K)\mathit{d}q$, 
              
              similiarly, the right hand side depends only on M, hence we can also consider it as a constant $C_{3}$.

                $\mathbf{First\; order\; estimate}$
                We want to give an estimate to $\mid\nabla v \mid$ by Schauder estimat,here L = $\Delta$ is an elliptic operator, and the equation Lv = $\Delta v$ is a uniformly elliptic equation(see G$\&$T), then by Schauder estimate we have
                
                $\mid\mid v\mid\mid _{C^{2,\alpha}}$  =  $\sum_{\mid \gamma\mid \le k}^{} \mid\mid D^{\gamma}v\mid\mid_{L^{\infty}} + \sum_{\mid \gamma\mid = k}^{}\mid D^{\gamma}v\mid_{\alpha}$
               
                $\le$ C($\mid\mid Lv\mid\mid_{C^{\alpha}}$ + $\mid\mid v\mid\mid_{L^{\infty}}$)
               
                =C($\mid\mid \Delta v\mid\mid_{L^{\infty}} + \mid\mid v\mid\mid_{L^{\infty}}$)
                
                 = C(sup$_{M\times [0,T)}$$\mid\Delta v\mid$ + sup$_{M\times [0,T)} \mid v\mid$ ), 
                
                while 
                
                sup$_{M\times [0,T)} \mid \nabla v\mid$ = $\mid\mid \nabla v\mid\mid_{L^{\infty}}$ $\le$
                $\sum_{\mid \gamma\mid \le k}^{} \mid\mid D^{\gamma}v\mid\mid_{L^{\infty}} + \sum_{\mid \gamma\mid = k}^{}\mid D^{\gamma}v\mid_{\alpha}$,  so
                
                sup$_{M\times [0,T)} \mid \nabla v\mid$ $\le$ C$_{6}$(sup$_{M\times [0,T)}$$\mid\Delta v\mid$ + sup$_{M\times [0,T)} \mid v\mid$)

                  where   $\mid\mid v\mid\mid _{C^{2,\alpha}}$  =  $\sum_{\mid \gamma\mid \le k}^{} \mid\mid D^{\gamma}v\mid\mid_{L^{\infty}} + \sum_{\mid \gamma\mid = k}^{}\mid D^{\gamma}v\mid_{\alpha}$  is the H$\ddot{o}$lder norms and  $\mid v\mid_{\alpha}$ is the $\alpha$-H$\ddot{o}$lder constant on $\Omega$, i.e.  $\mid f\mid_{\alpha}$  = sup$_{x\neq y \in \Omega}$ $\frac{\mid f(x)-f(y) \mid}{\mid x-y \mid^{\alpha}}$. 
                  
                therefore, we must give an estimate to 
                 $ sup_{M\times [0,T)}  \mid\Delta v \mid$  and $ sup_{M\times [0,T)} \mid v\mid $.

                 
                  $\mathbf{Lemma \;2}$

                  $\exists \; C_{1},C_{2}>0$, such that 0 $<$ n + $\Delta u$ $\le$ $C_{1} \mathit{e}^{C_{0}(u - inf_{M\times [0,T)}u)}$, for $\forall t \in [0,T)$

                 $ \mathit{Proof.}$ 
                 
                 Because $\tilde{\mathit{g_{i\bar{j}}}}$ - $\mathit{g_{i\bar{j}}}$ = $\frac{\partial^{2}u}{\partial\mathit{z^{i}}\partial\mathit{\bar{z}^{j}}}$, so we get 
                   $\mathit{g^{i\bar{j}}}$ $\tilde{\mathit{g_{i\bar{j}}}}$ = $\mathit{g^{i\bar{j}}}$($\mathit{g_{i\bar{j}}}$+$\frac{\partial^{2}u}{\partial\mathit{z^{i}}\partial\mathit{\bar{z}^{j}}}$)
                    = n + $\Delta u$, while $\mathit{g_{i\bar{j}}}$ and $\mathit{g^{i\bar{j}}}$ are both positive definite, so  $\mathit{g^{i\bar{j}}}$ $\tilde{\mathit{g_{i\bar{j}}}}$ is positive definite, 
                    
                    n +  $\Delta u$ is its trace so it's positive.

                    As for the other inequality, first we learn from Yau's work that:

                    If let $\tilde{\mathit{g_{i\bar{j}}}}$ =  $\mathit{g_{i\bar{j}}}$ + $\frac{\partial^{2} \varphi}{\partial\mathit{z^{i}}\partial\mathit{\bar{z}^{j}}}$ such that the equation holds
                    
                    det($\mathit{g_{i\bar{j}}}$ + $\frac{\partial^{2} \varphi}{\partial\mathit{z^{i}}\partial\mathit{\bar{z}^{j}}}$) $(det(\mathit{g_{i\bar{j}}}))^{-1}$ = $e^{F}$

                    then we get:

                    $\tilde{\Delta}$ ($e^{-C\varphi}$ (n + $\Delta \varphi$)) $\ge$ 
                    $e^{-C\varphi}$($\Delta F - n^{2}$ inf$_{i\ne l}$$(R_{i\tilde{i}l\tilde{l}}) $)
                    - C$e^{-C\varphi}$ n(n+$\Delta \varphi$) 
                    
                    +
                     (C + inf$_{i\ne l}$$(R_{i\tilde{i}l\tilde{l}}) $)$e^{-C\varphi}$ $e^{\frac{-F}{n-1}}(n+\Delta \varphi)^{\frac{n}{n-1}} $
                     
                     where $R_{i\tilde{i}l\tilde{l}} $ is the bisectional curvature of the  $\mathit{g_{i\bar{j}}}$, C$_{0}$ a positive constant such that
                     
                      C$_{0}$ + inf$_{i = 1}$ $(R_{i\tilde{i}l\tilde{l}}) $ $>$ 0.

                     Then we let 
                     
                     $\varphi$ = u (u satisfies the equation above), F = $\frac{\partial u}{\partial t}$-f, so we get:
                     
                     $\tilde{\Delta}$ ($e^{-Cu}$ (n + $\tilde{\Delta }u$)) $\ge$ 
                     
                      $e^{-C_{0}u}$($\Delta (\frac{\partial u}{\partial t} - f)$$ - n^{2} inf_{i\ne l}(R_{i\tilde{i}l\tilde{l}}) $) 
                     
                     - C$_{0}$ $e^{-C_{0}u}$n(n + $\Delta u$) 
                     
                     + (C$_{0}$ + inf$_{i\ne l}$$(R_{i\tilde{i}l\tilde{l}}) $) $e^{-C_{0}u}$ $e^{\frac{-(\frac{\partial u}{\partial t} - f)}{n-1}} $ (n+$\Delta u$)$^{\frac{n}{n-1}}$
                     
                     and  $\frac{\partial u}{\partial t}$ ($e^{-Cu}$ (n + $\tilde{\Delta }u$)) = $e^{-C_{0}u}$(-C$_{0}$$\frac{\partial u}{\partial t}$)(n + $\Delta u$) + $e^{-C_{0}u}$ $\Delta \frac{\partial u}{\partial t}$
                     
                     therefore, we get

                    $\diamond$ ($e^{-Cu}$ (n + $\tilde{\Delta }u$)) 
                    
                    $\ge$
                                       $e^{-C_{0}u}$($\Delta (\frac{\partial u}{\partial t} - f)$$ - n^{2} inf_{i\ne l}(R_{i\tilde{i}l\tilde{l}}) $) 
                    - C$_{0}$ $e^{-C_{0}u}$n(n + $\Delta u$)

                    + (C$_{0}$ + inf$_{i\ne l}$$(R_{i\tilde{i}l\tilde{l}}) $) $e^{-C_{0}u}$ $e^{\frac{-(\frac{\partial u}{\partial t} - f)}{n-1}} $ (n+$\Delta u$)$^{\frac{n}{n-1}}$
                    -$e^{-C_{0}u}$(-C$_{0}$$\frac{\partial u}{\partial t}$)(n + $\Delta u$) - $e^{-C_{0}u}$ $\Delta \frac{\partial u}{\partial t}$

                    =-$e^{-C_{0}u}$($\Delta f$ $ + n^{2} inf_{i\ne l}(R_{i\tilde{i}l\tilde{l}}) $) 
                    -C$_{0}$ $e^{-C_{0}u}$($n - \frac{\partial u}{\partial t}$)(n + $\Delta u$)

                    + (C$_{0}$ + inf$_{i\ne l}$$(R_{i\tilde{i}l\tilde{l}}) $) $e^{-C_{0}u}$ $e^{\frac{-\frac{\partial u}{\partial t} + f}{n-1}} $ (n+$\Delta u$)$^{\frac{n}{n-1}}$
                    
                    Then we assume for any t$\in$(0,T), the function ($e^{-Cu}$ (n + $\tilde{\Delta }u$))  achieves its maximum at (p$_{0}$,t$_{0}$)$\in$ M$\times$[o,t] and t$_{0}$$>$0, so at this point,
                    
                    $\diamond$ ($e^{-Cu}$ (n + $\tilde{\Delta }u$)) $\le$0
                     then 
                     
                     0$\ge$ -($\Delta f$ $ + n^{2} inf_{i\ne l}(R_{i\tilde{i}l\tilde{l}}) $) - C$_{0}$ ($n - \frac{\partial u}{\partial t}$)(n + $\Delta u$)

                     +(C$_{0}$ + inf$_{i\ne l}$$(R_{i\tilde{i}l\tilde{l}}) $) $e^{\frac{-\frac{\partial u}{\partial t} + f}{n-1}} $ (n+$\Delta u$)$^{\frac{n}{n-1}}$
                     
                     and from max$_{M}$$\mid$ $\frac{\partial u}{\partial t}$$\mid$ $\le$ max$_{M}$$\mid$f$\mid$ we get
                     
                     (n+$\Delta u$)$^{\frac{n}{n-1}}$ $\le$
                      $\frac{(\Delta f  + n^{2} inf_{i\ne l}(R_{i\tilde{i}l\tilde{l}}) ) + C_{0} (n - \frac{\partial u}{\partial t})(n + \Delta u)}{(C_{0} + inf_{i\ne l}(R_{i\tilde{i}l\tilde{l}}) ) e^{\frac{-\frac{\partial u}{\partial t} + f}{n-1}} }$

                    $\le$
                     $\frac{(\Delta f  + n^{2} inf_{i\ne l}(R_{i\tilde{i}l\tilde{l}}) ) + C_{0} (n + max_{M}\mid f\mid)(n + \Delta u)}{(C_{0} + inf_{i\ne l}(R_{i\tilde{i}l\tilde{l}}) )  }$,

                     while the$(\Delta f  + n^{2} inf_{i\ne l}(R_{i\tilde{i}l\tilde{l}}) )$,$C_{0} (n + max_{M}\mid f\mid)$$(C_{0} + inf_{i\ne l}(R_{i\tilde{i}l\tilde{l}}) )  $ are all independent of t, so there exists a constant C$^{'}$ such that

                   (n+$\Delta u$)$^{\frac{n}{n-1}}$ $\le$ C$^{'}+C^{'}(n + \Delta u)$
                   
                   then there exists a constant C$_{'1}$ such that
                   
                   $\frac{ (n+\Delta u)^{\frac{n}{n-1}}}{(1+n+\Delta u)}$ $\le$ C$_{'1}$,  which means $(n+\Delta u)$ $\le$ C$_{1}$ for a constant C$_{1}$ independent of t.
                   
                   Therefore, on M $\times$ [0,T) we have ($e^{-C_{0}u}$ (n + $\tilde{\Delta }u$))
                   
                    $\le$ $C_{1} e^{-C_{0} u(p,t_{0})}$ so n+$\Delta u$ 
                   
                   $<$ C$_{1} e^{C_{0}(u-u(p,t_{0}))}$

                   $\le$ C$_{1} e^{C_{0}(u-inf_{M\times [0,T)} u)}$  and C$_{0}$ and C$_{1}$ are both independent of t.


                    
                    $\mathbf{Lemma \;3}$
                    
                    There exists a constant $C_{4}$ so that $sup_{M\times [0,T)} \mid v \mid$ $<$ $C_{4}$.
                    
                    $ \mathit{Proof.}$
                    
                    Let 
                    $\omega$ = $\frac{\sqrt{-1}}{2}$ $\mathit{g_{i\bar{j}}}$ $\mathit{d}z^{i}$ $\wedge$  $\mathit{d}\bar{z}^{j}$,
                    $\tilde{\omega}$ = $\frac{\sqrt{-1}}{2}$ $\tilde{\mathit{g_{i\bar{j}}}}$ $\mathit{d}z^{i}$ $\wedge$  $\mathit{d}\bar{z}^{j}$, these are acutally the K$\ddot{a}$hler forms of $\mathit{g_{i\bar{j}}}$ and  $\tilde{\mathit{g_{i\bar{j}}}}$ ,separately. And the volume forms

                    $\mathit{d}$V = det($\mathit{g_{i\bar{j}}}$)  $\wedge^{n}_{i = 1}$ ( $\frac{\sqrt{-1}}{2}$ $\mathit{d}z^{i}$ $\wedge$  $\mathit{d}\bar{z}^{j}$) = $\frac{\omega ^n}{n!}$,

                    $\mathit{d}$$\tilde{V}$ = det( $\tilde{\mathit{g_{i\bar{j}}}}$)  $\wedge^{n}_{i = 1}$ ( $\frac{\sqrt{-1}}{2}$ $\mathit{d}z^{i}$ $\wedge$  $\mathit{d}\bar{z}^{j}$) = $\frac{\tilde{\omega} ^n}{n!}$
                    
                    because of the equation (3.1) we get
                    
                     log det($\tilde{\mathit{g_{i\bar{j}}}}$) - log det($\mathit{g_{i\bar{j}}}$ ) =  $\frac{\partial u}{\partial t}$ - f

                    then we get $\mathit{d}$$\tilde{V}$ = det($\mathit{g_{i\bar{j}}}$)$e^{\frac{\partial u}{\partial t} - f}$ $\wedge^{n}_{i = 1}$ ( $\frac{\sqrt{-1}}{2}$ $\mathit{d}z^{i}$ $\wedge$  $\mathit{d}\bar{z}^{j}$) = 
                    $e^{\frac{\partial u}{\partial t} - f}$ $\mathit{d}$V
                    
                    Therefore, for p$>$1,
                    
                    -$\frac{1}{n!}$$\int_{M}$$\frac{(-v)^{p-1}}{p-1}$ $(\omega ^n - \tilde{\omega} ^n )$ = -$\int_{M}$$\frac{(-v)^{p-1}}{p-1}$ $( \mathit{d}V - \mathit{d}\tilde{V} )$ 
                    
                    = $\int_{M}$$\frac{(-v)^{p-1}}{p-1}$ ($e^{\frac{\partial u}{\partial t} - f} - 1$)   $\mathit{d}$V

                    Because from lemma 1 we know v is bounded then we can renormalized it so that v$<$-1.
                    Then on the other hand,

                    -$\int_{M}$$\frac{(-v)^{p-1}}{p-1}$ $(\omega ^n - \tilde{\omega} ^n )$

                     =$\int_{M}$$\frac{(-v)^{p-1}}{p-1}$ ($\tilde{\omega}$ - $\omega$)$\sum_{j = 0}^{ n-1}$ $\tilde{\omega}^{j}$ $\wedge$$\omega^{n-j-1}$

                      =  $\int_{M}$$\frac{(-v)^{p-1}}{p-1}$ ($\frac{\sqrt{-1}}{2}$$\partial \bar{\partial}v$)$\wedge$$\sum_{j = 0}^{ n-1}$ $\tilde{\omega}^{j}$ $\wedge$$\omega^{n-j-1}$
                      
                      = $\int_{M}$ $\frac{(-v)^{p-1}}{p-1}$   ($\frac{\sqrt{-1}}{2}$$\mathit{d\frac{1}{2}(\partial v + \bar{\partial}v)}$)  $\wedge$$\sum_{j = 0}^{ n-1}$ $\tilde{\omega}^{j}$ $\wedge$$\omega^{n-j-1}$
                      
                      =$\int_{M}$ $\mathit{d}$$\frac{(-v)^{p-1}}{p-1}$   ($\frac{\sqrt{-1}}{2}$$\frac{1}{2}(\partial v + \bar{\partial}v)$)  $\wedge$$\sum_{j = 0}^{ n-1}$ $\tilde{\omega}^{j}$ $\wedge$$\omega^{n-j-1}$

                      =$\int_{M}$ $(-v)^{p-2}$ $\mathit{d}v$  ($\frac{\sqrt{-1}}{2}$$\frac{1}{2}(\partial v + \bar{\partial}v)$)  $\wedge$$\sum_{j = 0}^{ n-1}$ $\tilde{\omega}^{j}$ $\wedge$$\omega^{n-j-1}$
                      
                      =$\int_{M}$ $(-v)^{p-2}$  ($\frac{\sqrt{-1}}{2}$$(\partial v \wedge \bar{\partial}v)$)  $\wedge$$\sum_{j = 0}^{ n-1}$ $\tilde{\omega}^{j}$ $\wedge$$\omega^{n-j-1}$
                      
                      there we integral by part because v vanishes on the boundry of M.
                      
                      then the integral above
                      
                      $\ge$ 
                      $\int_{M}$ $(-v)^{p-2}$  ($\frac{\sqrt{-1}}{2}$$(\partial v \wedge \bar{\partial}v)$)  $\wedge$ $\omega ^{n-1}$
                      because each term in  
                      
                       $\frac{\sqrt{-1}}{2}$$(\partial v \wedge \bar{\partial}v)$  $\wedge$ $\tilde{\omega}^{j}$ $\wedge$$\omega^{n-j-1}$ is nonnegative.
                      
                      Then let $\mid \nabla v\mid^{2}$ =$\mathit{g^{i\bar{j}}}$ $\frac{\partial v}{\partial z^{i}}$ $\frac{\partial v}{\bar{\partial z}^{i}}$  we have
                      
                       $\int_{M}$ $(-v)^{p-2}$ $\mid \nabla v\mid^{2}$ $\mathit{d}V$ $\le$
                       n! $\int_{M}$$\frac{(-v)^{p-1}}{p-1}$ ($e^{\frac{\partial u}{\partial t} - f} - 1$)   $\mathit{d}$V.
                       
                       While because (-v$^{p-2}$)$\mid \nabla v\mid$$^{2}$ =  4p$^{-2}$ $\mid \nabla(-v)^{\frac{p}{2}}\mid$$^{2}$
                       we replace the corresponding term in the inequality above:
                       
                       $\int_{M}$ 4p$^{-2}$$\mid \nabla(-v)^{\frac{p}{2}}\mid$$^{2}$$\mathit{d}V$ $\le$
                        n! $\int_{M}$$\frac{(-v)^{p-1}}{p-1}$ ($e^{\frac{\partial u}{\partial t} - f} - 1$)   $\mathit{d}$V, then
                        
                        $\int_{M}$ $\mid \nabla(-v)^{\frac{p}{2}}\mid$$^{2}$ $\mathit{d}V$ $\le$ $\frac{n!}{4}$ p$^{2}$$\int_{M}$ $\frac{(-v)^{p-1}}{p-1}$ ($e^{\frac{\partial u}{\partial t} - f} - 1$)   $\mathit{d}$V
                        
                        while due to the compability condition, we get the term ($e^{\frac{\partial u}{\partial t} - f} - 1$) is positive and bounded, therefore, the term above :
                        
                        $\le$ C$\frac{p^{2}}{p-1}$ $\int_{M}$ $(-v)^{p-1}$ $\mathit{d}V$
                        
                        and by the norm $\mid\mid$ $(-v)^{\frac{p}{2}}$$\mid\mid$$^{2}$$_{H^{1}}$ 
                        
                        = 
                         $\int_{M}$ $\mid \nabla(-v)^{\frac{p}{2}}\mid$$^{2}$ $\mathit{d}V$ + $\int_{M}$$(-v)^{p}$ $\mathit{d}V$

                         $\le$ (C$\frac{p^{2}}{p-1}$ + 1) $\int_{M}(-v)^{p}$$\mathit{d}V$

                         $\le$ Cp $\int_{M}(-v)^{p}$$\mathit{d}V$ , for p$>$ 1, there the C has changed but we still use it for a constant. When p = 1 then we just replace $\frac{(-v)^{p-1}}{p-1}$ by the term log(-v) and the process works well.
                         
                          -$\int_{M}$ $log(-v)$ $(\omega ^n - \tilde{\omega} ^n )$ =
                         $\int_{M}$$\frac{(-v)^{p-1}}{p-1}$ ($\tilde{\omega}$ - $\omega$)$\sum_{j = 0}^{ n-1}$ $\tilde{\omega}^{j}$ $\wedge$$\omega^{n-j-1}$

                         =  $\int_{M}$$log(-v)$ ($\frac{\sqrt{-1}}{2}$$\partial \bar{\partial}v$)$\wedge$$\sum_{j = 0}^{ n-1}$ $\tilde{\omega}^{j}$ $\wedge$$\omega^{n-j-1}$
                         
                         = $\int_{M}$ $log(-v)$   ($\frac{\sqrt{-1}}{2}$$\mathit{d\frac{1}{2}(\partial v + \bar{\partial}v)}$)  $\wedge$$\sum_{j = 0}^{ n-1}$ $\tilde{\omega}^{j}$ $\wedge$$\omega^{n-j-1}$
                         
                         =$\int_{M}$ $\mathit{d}$$log(-v)$   ($\frac{\sqrt{-1}}{2}$$\frac{1}{2}(\partial v + \bar{\partial}v)$)  $\wedge$$\sum_{j = 0}^{ n-1}$ $\tilde{\omega}^{j}$ $\wedge$$\omega^{n-j-1}$

                         =-$\int_{M}$ $(-v)^{-1}$ $\mathit{d}v$  ($\frac{\sqrt{-1}}{2}$$\frac{1}{2}(\partial v + \bar{\partial}v)$)  $\wedge$$\sum_{j = 0}^{ n-1}$ $\tilde{\omega}^{j}$ $\wedge$$\omega^{n-j-1}$
                         
                         =-$\int_{M}$ $(-v)^{-1}$  ($\frac{\sqrt{-1}}{2}$$(\partial v \wedge \bar{\partial}v)$)  $\wedge$$\sum_{j = 0}^{ n-1}$ $\tilde{\omega}^{j}$ $\wedge$$\omega^{n-j-1}$
                         
                         While -$\frac{1}{n!}$$\int_{M}$$log(-v)$ $(\omega ^n - \tilde{\omega} ^n )$ 
                         
                         = -$\int_{M}$$log(-v)$ $( \mathit{d}V - \mathit{d}\tilde{V} )$ 
                         
                         = $\int_{M}$$log(-v)$ ($e^{\frac{\partial u}{\partial t} - f} - 1$)   $\mathit{d}$V
                         
                         So   $\int_{M}$ $(-v)^{p-1}$ $\mid \nabla v\mid^{2}$ $\mathit{d}V$ $\le$
                         n! $\int_{M}$ $log(-v)$ ($e^{\frac{\partial u}{\partial t} - f} - 1$)   $\mathit{d}$V.
                         
                        and  (-v$^{-1}$)$\mid \nabla v\mid$$^{2}$ =  4 $\mid \nabla(-v)^{\frac{1}{2}}\mid$$^{2}$
                         
                         Then  $\mid\mid$ $(-v)^{\frac{1}{2}}$$\mid\mid$$^{2}$$_{H^{1}}$ 
                         
                         = 
                         $\int_{M}$ $\mid \nabla(-v)^{\frac{1}{2}}\mid$$^{2}$ $\mathit{d}V$ + $\int_{M}$$(-v)^{1}$ $\mathit{d}V$

                         $\le$ C$\int_{M} log(-v)$ $\mathit{d}V$ + $\int_{M}(-v)^{p}$$\mathit{d}V$

                         $\le$ Cp $\int_{M}(-v)^{p}$ $\mathit{d}V$ for v$<$-1 we assumed before and log(-v) is slower than(-v).
                        So the equality works well when p$\ge$1.

                         Because $\int_{M}v\mathit{d}V$ = 0,  we learn from the Sobolov inequality(G$\&$T[5] p155)
                         
                         $\mid\mid(-v)^{\frac{p}{2}}\mid\mid _{L^{\frac{n}{n-1}}}$ 
                         
                         $\le$ C$\mid\mid D((-v)^{\frac{p}{2}})\mid\mid_{H^{1}}$ 
                         
                         = C$\mid\mid \frac{p}{2} (-v)^{-\frac{p}{2}} (-Dv)\mid\mid_{H^{1}}$ while as our assume v<-1 and we have an estimate for $\mid\mid \nabla v\mid\mid$ so we get $\mid\mid(-v)^{-\frac{p}{2}}\mid\mid$ $\le$ $\mid\mid(-v)^{\frac{p}{2}}\mid\mid$   and then

                         C$\mid\mid \frac{p}{2} (-v)^{-\frac{p}{2}} (-Dv)\mid\mid_{H^{1}}$ $\le$ C$\mid\mid (-v)^{\frac{p}{2}}\mid\mid_{H^{1}}$, therefore,
                         
                         $\mid\mid(-v)^{\frac{p}{2}}\mid\mid^{2} _{L^{\frac{2n}{n-1}}}$ $\le$  $\mid\mid(-v)^{\frac{p}{2}}\mid\mid^{2} _{L^{\frac{2n}{n-1}}}$
                         $\le$ C$\mid\mid (-v)^{\frac{p}{2}}\mid\mid^{2}_{H^{1}}$

                         While  $\int_{M}(-v)^{p}$$\mathit{d}V$ gives us the L$^{p}$ norm,so we get 
                         
                         $\mid\mid v\mid\mid^{p}$${_{L}}^{\frac{np}{n-1}}$ $\le$ Cp $\mid\mid v\mid\mid ^{p}$ ${_{L}}^{p}$ for any p$\ge$1
                         
                         Then we consider $\gamma$ = $\frac{n}{n-1}$ and p = $\gamma^{j}$ for j = 0,1,2,......,
                         then 
                         
                         $\mid\mid v\mid\mid$ ${^{\gamma}}^{j} {_{L}}^{\gamma^{j+1}}$ $\le$ $C\gamma^{j} $
                         $\mid\mid v\mid\mid$ ${^{\gamma}}^{j} {_{L}}^{\gamma^{j}}$,
                         we continue decrease j so:
                         
                         $\mid\mid v\mid\mid$ $ {_{L}}^{\gamma^{j+1}}$ $\le$ $C^{-\gamma^{j}}\gamma^{\frac{j}{\gamma^{j}}} $
                         $\mid\mid v\mid\mid$ $ {_{L}}^{\gamma^{j}}$ ;
                         
                         $\mid\mid v\mid\mid$ $ {_{L}}^{\gamma^{j+1}}$
                          $\le$ $C^{-\gamma^{j}-\gamma^{j-1}}\gamma^{\frac{j}{\gamma^{j}}+\frac{j-1}{\gamma^{j-1}}} $
                         $\mid\mid v\mid\mid$ $ {_{L}}^{\gamma^{j-1}}$ ;
                          $\mid\mid v\mid\mid$ $ {_{L}}^{\gamma^{j+1}}$$\le$...;

                         this process continue until 0, i.e.
                         
                          $\mid\mid v\mid\mid$ $ {_{L}}^{\gamma^{j+1}}$ $<$
                          C$^{\sum_{k =0}^{j}\frac{1}{\gamma^{k}}}$ $\gamma^{\sum_{k =0}^{j}\frac{k}{\gamma^{k}}}$ C$_{3}$, for L$_{0}$ norm is a constant.
                         
                         Finally, we let j increases to infinity then $\mid\mid v\mid\mid$ ${_{L}}^{\infty}$ $\le$ $C_{4}$ for C$_{4}$ a constant from the right hand side of the inequality induced above so it's independent of t:
                         
                         sup$_{M\times [0,T)}$ $\mid v\mid$ $\le$ C$_{4}$, this finishes the proof.

                         After the long and complex calculation, we can continue our estimate.

                        Because $\frac{1}{Vol(M)}$$\int_{M}u\mathit{d}V$ is not relative to $\partial \tilde{\partial}$, so we get:

                         0$<$n+$\Delta v$ = n+$\Delta u$ $\le$ $C_{1} \mathit{e}^{C_{0}(u - inf_{M\times [0,T)}u)}$
                         = $C_{1} \mathit{e}^{C_{0}(v - inf_{M\times [0,T)}v)}$
                         $\le$ C$_{5}$

                         so the $\Delta v$ is bounded, the equality holds because the difference is not relevant to the constant $\frac{1}{Vol(M)}$$\int_{M}u\mathit{d}V$, then the constant C$_{5}$ appears because sup$\mid v\mid$ is bounded.
                         
                         The by the Schauder estimate we presented before

                         sup$_{M\times[0,T)}$ $\mid \nabla v\mid$ $\le$
                         C$_{6}$(sup$_{M\times[0,T)}$ $\mid \Delta v\mid$ + sup$_{M\times[0,T)}$ $\mid v\mid$) $\le$
                         C$_{6}$(constant +constant) by Lemma 2 and Lemma 3 then
                         sup$_{M\times[0,T)}$ $\mid \nabla v\mid$ $\le$
                         C$_{7}$. This finishes the first order estimate.

                         $\mathbf{Second\; order\; estimate }$

                         From the estimate of  n+$\Delta u $, in Yau's work the choice of the metric(in page 348) on one hand make the matrix ($\delta $ +u$_{i\tilde{i}}$) is positive definite and Hermitian so we know that 1+u$_{i\tilde{i}}$ is bounded above for any i, by the metric Yau chose, $\prod_{i=1}^{m} (1+u_{i\tilde{i}})$ = $e^{\frac{\partial u}{\partial t} - f}$ gives a lower estimate which is positive, this is because the product of these terms has a upper bound, so each term can not go to the negative infinity, so they all have a lower bound, while the matrix ($\delta $ +u$_{i\tilde{i}}$) is positive definite and Hermitian so each term is positive,  therefore,  there exists two positive constant A, B such that 
                         
                         A$\le$1+u$_{i\tilde{i}}$ $\le$B for any i.

                          $\mathbf{Third\; order\; estimate }$

                          First let S = $\sum_{}^{}$ ${\tilde{g}}^{i\bar{r}}$
                          ${\tilde{g}}^{j\bar{s}}$ ${\tilde{g}}^{k\bar{t}}$
                          v$_{i\bar{j}k}$ v$_{\bar{r}s\bar{t}}$. Actually in this disturbing definition, Cao said in his paper that he followed E.Calabi and Yau, Yau said he followed E.Calabi, but I have not read E.Calabi's paper, I have only read Yau and Cao's paper, so I do not know who I am following, that's really insteresting. To show respect to E.Calabi, I decide to write I'm following E.Calabi here, too. I'll read his work in a few minutes.
                          
                          We note A$\simeq $B if $\mid A-B\mid$ $\le$ C$_{1}$ $\sqrt{S}$+$C_{2}$ for C$_{1}$ and C$_{2}$ are constants can be estimated, further more, we note A$\cong$B if $\mid A-B\mid$ $\le$ C$_{3}$S + C$_{4}$ $\sqrt{S}$+$C_{5}$ for C$_{3}$ and C$_{4}$, C$_{5}$ are constants can be estimated. And again by the metric Yau defined before, through totally 60 rows calculations in Yau's appendix, 
                          
                          $\tilde{\Delta}$S $\simeq$ $\sum_{}^{}$ (1+v$_{i\bar{i}}$)$^{-1}$ (1+v$_{j\bar{j}}$)$^{-1}$
                          (1+v$_{k\bar{k}}$)$^{-1}$ (1+v$_{\alpha\bar{\alpha}}$)$^{-1}$ $\times$
                          \{ $\mid$   v$_{\bar{i}j\bar{k}\alpha}$ - $\sum_{}^{}$ v$_{\bar{i}b\bar{k}}$ v$_{\bar{p}j\alpha}$ (1+v$_{p\bar{p}}$)$^{-1}$ $\mid$$^{2}$ +$\mid$ v$_{i\bar{j}k\alpha}$ -$\sum_{p}^{}$ (v$_{\bar{p}i\alpha}$  v$_{p\bar{j}k}$ + v$_{\bar{p}ik}$  v$_{p\bar{j}\alpha}$)(1+v$_{p\bar{p}}$)$^{-1}$  $\mid$ $^{2}$      \}.
                          
                          And $\tilde{\Delta}$($\Delta v$) $\ge$ $\sum_{}^{}$ (1+v$_{k\bar{k}}$)$^{-1}$ (1+v$_{i\bar{i}}$)$^{-1}$
                          $\mid v_{k\bar{ij}}\mid$$^{2}$ - C$_{6}$, where C$_{6}$ is  a constant which can be estimated, then without loss generality, we can choose a big C$_{7}$ such that
                          
                           $\tilde{\Delta}$(S+C$_{7}$ $\Delta$v) $\ge$
                          C$_{8}$S - C$_{9}$, for C$_{7}$,C$_{8}$,C$_{9}$, are positive constants can be estimated.

                         Then we assume p(t) is the maximum point of the funciton S+C$_{7}$ $\Delta v$, we get
                         
                         0$\ge$C$_{8}$ S-C$_{9}$; so  C$_{8}$ S$\le$C$_{9}$;
                         
                         C$_{8}$(S+C$_{7}\Delta v$) $\le$ C$_{9}$ + C$_{8}$C$_{7}$$\Delta v$.
                         
                         Due to the $\Delta v $ is bounded from estimate before, this gives an estimate to sup$_{M\times[0,T)}$(M+C$_{7}$$\Delta v$) therefore, we have an estimate of sup$_{M}$S,as Yau said, this gives an estimate to v$_{i\bar{j}k}$ and similar term in terms of g$_{i\bar{j}}$, sup$\mid F\mid$,
                         
                         sup$\mid \nabla F\mid$,sup$_{M}$sup$_{i}$$\mid F_{i\bar{i}}\mid$ and sup$_{M}$sup$_{i,j,k}$$\mid F_{i\bar{j}k}\mid$ .
                         
                      $   \mathbf{Long\; time\;existence}$

                         Finally we prove the proposition.
                         
                         $\mathbf{Proposition\; of \; the \; long \; existence}$
                         
                         Assume u be the solution of the equation

                         $\frac{\partial u}{\partial t}$ =
                         log det($\mathit{g_{i\bar{j}}}$ + $\frac{\partial^{2}u}{\partial\mathit{z^{i}}\partial\mathit{\bar{z}^{j}}}$ ) - log det($\mathit{g_{i\bar{j}}}$ ) + f----------(3.1)
                         
                         where t$\in$ [0,T) which is the maximum time interval. Let v be the normalization of u :

                         v = u - $\frac{1}{Vol(M)}$$\int_{M}u\mathit{d}V$. Then the C$_{\infty}$ norm of v are uniformly bounded for any t$\in$(0,T), which means T = $\infty$.
                         Then there exists a sequence t$_{n}$ increasing to infinity such that v(x,t$_{n}$) converges in the topology generated by C$_{\infty}$ norm to a smooth function v$_{\infty}$(x) on M when n increases to infinity.
                         
                         $\mathit{Proof.}$
                         
                         Differentiate the equation(3.1):
                         
                         $\frac{\partial}{\partial t}$ ($\frac{\partial u}{\partial z^{k}}$) = $\tilde{\mathit{g}}^{i\bar{j}}$($\frac{\partial}{\partial z^{k}}$$\mathit{g_{i\bar{j}}}$ + $\partial \bar{\partial}$($\frac{\partial u}{\partial z^{k}}$)) - $\mathit{g^{i\bar{j}}}$ $\frac{\partial}{\partial z^{k}}$ $\mathit{g_{i\bar{j}}}$ + $\frac{\partial f}{\partial z^{k}}$;
                         
                         so $\diamond $ ($\frac{\partial u}{\partial z^{k}}$)= $\mathit{g^{i\bar{j}}}$ $\frac{\partial}{\partial z^{k}}$ $\mathit{g_{i\bar{j}}}$ - $\tilde{\mathit{g}}^{i\bar{j}}$$\frac{\partial}{\partial z^{k}}$$\mathit{g_{i\bar{j}}}$ + $\frac{\partial f}{\partial z^{k}}$;
                         
                         we have the estimate of $\frac{\partial u}{\partial t}$ and $\tilde{\Delta u}$ so the coeffecients of the $\diamond$ are bounded, and by the C$^{0,\alpha}$ norm is defined as: $\mid\mid u \mid\mid_{C^{0,\alpha}}$ = $\mid u\mid_{L^{\infty}}$ + $\mid u\mid_{\alpha}$, according to our 0 order estimate, they are bounded in the $C^{0,\alpha}$ norm; and these term are actually the H$\ddot{o}$lder coefficient so is has bounded C$^{0;\alpha}$ norm. While the right hand side have similar estimate for C$^{0;\alpha}$ norm for 0$<\alpha<1$ because the RHS only depends on the Kahler metric and f which are smooth.
                         By interior regularity theory(see G$\&$ T section 6.4) which gives the C$^{k+2,\alpha}$ estimate when then coeffecients and the nonhomogeneous term are of C$^{k,\alpha}$ estimate, this $\diamond$ is elliptic and the RHS terms are smooth so they have C$^{0,\alpha}$ estimate, then as the solution of this equation, $\frac{\partial u}{\partial z^{k}}$ has the C$^{2}$ estimate we proved before, then $\frac{\partial u}{\partial z^{k}}$ has a uniform C$^{2,\alpha}$ estimate, similarly,$\frac{\partial u}{\bar{\partial} \bar{z}^{k}}$ also has this estimate,  there the C$^{0;\alpha}$ norm estimate gives the C$^{2;\alpha}$ norm estimate, so $\frac{\partial u}{\partial z^{k}}$ and $\frac{\partial u}{\partial \bar{z}^{k}}$ are uniformly bounded in C$^{2;\alpha}$ norm estimate, then because 
                             $\tilde{\mathit{g_{i\bar{j}}}}$ = $\mathit{g_{i\bar{j}}}$ + $\frac{\partial^{2}u}{\partial\mathit{z^{i}}\partial\mathit{\bar{z}^{j}}}$
                             then the C$^{2,\alpha}$ estimate of $\frac{\partial u}{\partial z^{k}}$ actually gives the C$^{1,\alpha}$ estimate of the RHS and the coeffecients of $\diamond$ which are determined by the metrics, therefore, the right hand side and the coefecients of $\diamond$ have uniform C$^{1,\alpha}$ estimate .
                         Use the Interior regularity theory again, then so $\frac{\partial u}{\partial z^{k}}$ and $\frac{\partial u}{\partial \bar{z}^{k}}$ are uniformly bounded in C$^{3;\alpha}$ norm estimate, then again by our fundamental equation of u and metrics , the coeffecients and the RHS have anoter one order more estimate, then $\frac{\partial u}{\partial z^{k}}$ and $\frac{\partial u}{\partial \bar{z}^{k}}$ have  another two orders more estimate......  Then we repeat use the theory, by iteration, then v(x,t) has uniformly bounded C$^{\infty}$ norm for any t $\in$ (0,T), finally we choose a sequence of t goes into the infinity such that it has a subsequence t$_{n}$ making v(x,t$_{n}$) converge to a smooth function v$_{\infty}$ , so that solution exists. And then because $\frac{\partial u}{\partial t}$ is uniformly bounded referring to t, and as t goes into the infinity u can not blow up in finite time, this means that our estimates are independent of t then when solution metric get the upper bound of t, says T, then our estimates still works the estimates above are independent of t so if we choose t$_{0}\in[0,T)$ then the solution also exists in [t$_{0}$,t$_{0}+\epsilon$] for $\epsilon$ independent of t, so we can use this  $\tilde{\mathit{g_{i\bar{j}}}}$(T) as the intial condition of the same equation chosing a new initial point and continue to deformation, so the process can continue to infinity because our estimates always work,  so u exists for all time.
                         
                         Up to here, the long time existence has finally been proved.

                         \section{Arguement for uniform convergence}

                         At this time, since we have proved the existence of the long time solution, then to follow our idea, we need to show the uniform convergence of the solution and then let the time goes into the infinity to get the final conclusion.

                         $\frac{\partial u}{\partial t}$ =
                         log det($\mathit{g_{i\bar{j}}}$ + $\frac{\partial^{2}u}{\partial\mathit{z^{i}}\partial\mathit{\bar{z}^{j}}}$ ) - log det($\mathit{g_{i\bar{j}}}$ ) + f--------------------(3.1)
                         
                         and initial status: u(x,t) = 0 when t = 0 on M$\times$[0,$\infty$), similarly, we renormalized u:

                          v = u - $\frac{1}{Vol(M)}$$\int_{M}u\mathit{d}V$ 
                          
                          And then we'll show the uniform convergence of v(x,t) and $\frac{\partial u}{\partial t}$ when t goes into the infinity.

                         We learn from before that  
                         $\frac{\partial}{\partial t}$($\frac{\partial u}{\partial t}$) = ${\tilde\mathit{g}^{i\bar{j}}}$ $\frac{\partial^{2}}{\partial\mathit{z^{i}}\partial\mathit{\bar{z}^{j}}}$($\frac{\partial u}{\partial t}$) so ($\tilde{\Delta}$ - $\frac{\partial}{\partial t}$) $\frac{\partial u}{\partial t}$ = 0 
                         
                         where  $\frac{\partial u}{\partial t}$ (x,t) =f(x) when t = 0. So in order to analyze u, we should first take care of this equation 
                         
                         ($\tilde{\Delta}$ - $\frac{\partial}{\partial t}$) $\frac{\partial u}{\partial t}$ = 0
                         
                         And following Yau's work, Cao gave a modification of an important theory of this equation.

                       $\mathbf{Theorm\;4.1}$
                       
                       We assume M be a compact manifold whose dimension is n, $\mathit{g_{i\bar{j}}}$(t)  be a family of Riemannian metrics on M, such that the following holds:
                       
                       (1).$\exists$ constants C$_{1}$,C$_{2}$ positive and independent of t such that 
                       
                       C$_{1}$  $\mathit{g_{ij}}$(0) $\le$  $\mathit{g_{ij}}$(t) $\le$  C$_{2}$ $\mathit{g_{ij}}$(0)
                       
                       (2).$\exists$ constants C$_{3}$ positive and independent of t such that
                       
                       $\mid \frac{\partial \mathit{g_{ij}} }{\partial t} \mid$(t) $\le$ C$_{3}$ $\mathit{g_{ij}}$(0)
                       
                       (3).$\exists$ constants K positive and independent of t such that
                       
                       R$_{ij}$(t) $\ge$ -K  $\mathit{g_{ij}}$(0)
                       
                       Then we assume $\phi$ is positive and satisfies the equation:
                       
                        ($\Delta_{t}$ - $\frac{\partial}{\partial t}$) $\frac{\partial \phi}{\partial t}$ = 0
                        
                        on M$\times[0,\infty)$ where $\Delta_{t}$ is the Laplace operator, then  $\forall \alpha$ $>$1,
                        
                        sup$_{x\in M}$ $\phi$(x,t$_{1}$) $\le$ inf$_{x\in M}$ $\phi$ (x,t$_{2}$) ($\frac{t_{2}}{t_{1}}$)$^{\frac{n}{2}}$ e$^{(
                        	\frac{1}{4(t_{2} - t_{1})} C_{2}^{2} \mathit{d}^{2} 
                        	+ (\frac{n\alpha K}{2(\alpha-1)} + C_{2}C_{3}(n + A) )(t_{2} - t_{1})
                        	)}$
                        
                        for $\mathit{d}$ is the diameter of M measured by   $\mathit{g_{ij}}$(0), i.e. 
                        $\mathit{d}$ = sup$_{x,y\in M}$ g$_{ij}$(0)(x,y);  and 
                        
                        A = sup$\mid \mid \nabla^{2} log \phi \mid \mid$; and 0$< t_{1}< t_{2}<\infty$.
                        
                        Cao did not show the proof because the proof of the theorm is totally a tough work in Yau's paper. Now, we can use this conclusion, let F = $\frac{\partial u}{\partial t}$, then by the maximum principle we still consider t=0 as the boundry of M$\times$[0,$\infty$] for this parabolic equation and $t_{2}$ $> t_{1} >0$, we get:
                        
                        sup$_{x\in M}$ F(x,t$_{2}$) $<$  sup$_{x\in M}$ F(x,t$_{1}$) $<$  sup$_{x\in M}$ f(x)

                       inf$_{x\in M}$ F(x,t$_{2}$) $>$  inf$_{x\in M}$ F(x,t$_{1}$) $>$  inf$_{x\in M}$ f(x)
                       
                       here $t_{2}$ $> t_{1} >0$ because we can always choose at $t_{2}$ the $\partial u$ converges more than $t_{1}$.
                       
                       $\mathbf{Remark}$
                       The conditions above also hold for $\tilde{\mathit{R_{i\bar{j}}}}$, i.e.

                       (1).$\exists$ constants C$_{1}$,C$_{2}$ positive and independent of t such that 
                       
                       C$_{1}$   $\tilde{\mathit{g_{i\bar{j}}}}$(0) $\le$   $\tilde{\mathit{g_{i\bar{j}}}}$(t) $\le$  C$_{2}$  $\tilde{\mathit{g_{i\bar{j}}}}$(0)
                       
                       (2).$\exists$ constants C$_{3}$ positive and independent of t such that
                       
                       $\mid \frac{\partial  \tilde{\mathit{g_{i\bar{j}}}} }{\partial t} \mid$(t) $\le$ C$_{3}$  $\tilde{\mathit{g_{i\bar{j}}}}$(0)
                       
                       (3).$\exists$ constants K positive and independent of t such that
                       
                   $\tilde{\mathit{R_{i\bar{j}}}}$(t) $\ge$ -K   $\tilde{\mathit{g_{i\bar{j}}}}$(0)

                        Actually,
                           $\tilde{\mathit{g_{i\bar{j}}}}$(0) = ($\mathit{g_{ij}}$(0) + $\frac{\partial^{2}u}{\partial\mathit{z^{i}}\partial\mathit{\bar{z}^{j}}}$(0) ),
                         from the estimate before we know $\frac{\partial^{2}u}{\partial\mathit{z^{i}}\partial\mathit{\bar{z}^{j}}}$(0) is bounded,then
                         we can always choose a  C$_{1}$$\le$ $\frac{ (\mathit{g_{ij}}(t) + \frac{\partial^{2}u}{\partial\mathit{z^{i}}\partial\mathit{\bar{z}^{j}}}(t) )}{ (\mathit{g_{ij}}(0) + \frac{\partial^{2}u}{\partial\mathit{z^{i}}\partial\mathit{\bar{z}^{j}}}(0) )}$   because 
                         $\mathit{g_{ij}}$(0) $\le$ C$\tilde{\mathit{g_{i\bar{j}}}}$(c) for some C, similary, choose
                         
                          C$_{2}$$\ge$ $\frac{ (\mathit{g_{ij}}(t) + \frac{\partial^{2}u}{\partial\mathit{z^{i}}\partial\mathit{\bar{z}^{j}}}(t) )}{ (\mathit{g_{ij}}(0) + \frac{\partial^{2}u}{\partial\mathit{z^{i}}\partial\mathit{\bar{z}^{j}}}(0) )}$ , so the (1) holds;

                          then  $\mid \frac{\partial  \tilde{\mathit{g_{i\bar{j}}}} }{\partial t} \mid$(t) $\le$  $\mid \frac{\partial \mathit{g_{ij}} }{\partial t} \mid$(t) + $\mid \partial\tilde{\partial}(\frac{\partial u}{\partial t})\mid$, while the estimate before tells that
                          $\mid \partial\tilde{\partial}(\frac{\partial u}{\partial t})\mid$(t) and $\partial\tilde{\partial} u$(0) are bounded
                          so we can choose a C$_{3}$ such that 
                          
                          $\partial\tilde{\partial} u$(0) $\ge$ C$_{3}$ $\mid \partial\tilde{\partial}(\frac{\partial u}{\partial t})\mid$(t), so (2) holds;
                          (3) Vy the long time existence thoerem,
                          
                           $\tilde{\mathit{g}_{ij}}(0) $ = $\mathit{g_{ij}} $ so when it holds for $\mathit{g_{ij}} $ then in the convergence process, $\tilde{\mathit{g}_{ij}}(t) $ converges to $\tilde{\mathit{g}_{ij}}(\infty) $, so whenever t$\in$[0,$\infty$], there exists such K.

                          Then we define 
                          
                          $\varphi_{n}$(x,t) = sup$_{x\in M}$F(x,n-1) - F(x,n-1+t)
                          
                          $\phi_{n}$(x,t) = F(x,n-1+t) - inf$_{x\in M}$F(x,n-1)
                          
                          $\omega$(t) = sup$_{x\in M}$F(x,t) - inf$_{x\in M}$F(x,t)
                          
                          While because the sup and the inf are constants,and the $\mathit{d}$(n-1)+t =$\mathit{d}$t so the $\varphi_{n}$(x,t) and $\phi_{n}$(x,t) both satisfy the equation and the initial condition, and by the ineqaulity above, they are both positive. Then we take t$_{1}$ = $\frac{1}{2}$, t$_{2}$ = 1, then we use the Theorm 4.1 in  $\varphi_{n}$(x,t) and $\phi_{n}$(x,t) separately:
                          
                           sup$_{x\in M}$ $\phi_{n}$(x,$\frac{1}{2}$) $\le$  inf$_{x\in M}$ $\phi_{n}$ (x,1)$\gamma$, 
                           
                           where $\gamma$ = 2$^{\frac{n}{2}}$ e$^{\frac{1}{2} C_{2}^{2} \mathit{d}^{2}  + \frac{1}{2} (\frac{n\alpha K}{2(\alpha-1)} + C_{2}C_{3}(n+A) )}$ is independent of t, then
                           
                           sup$_{x\in M}$F(x,n-1)  - inf$_{x\in M}$F(x,n-$\frac{1}{2}$) $\le$ $\gamma$(sup$_{x\in M}$F(x,n-1)  - sup$_{x\in M}$F(x,n))
                           
                           similarly,
                           
                           sup$_{x\in M}$F(x,n$\frac{1}{2}$)  - inf$_{x\in M}$F(x,n-1) $\le$ $\gamma$(inf$_{x\in M}$F(x,n)  - inf$_{x\in M}$F(x,n-1))
                           
                           do not forget we have $\omega$(t) = sup$_{x\in M}$F(x,t) - inf$_{x\in M}$F(x,t) so we add the two inequality together:
                           
                           $\omega(n-1) + \omega(n-\frac{1}{2})$ $\le$ $\gamma$($\omega(n-1) - \omega(n)$)
                           
                           because $\omega(n)$ $\ge$0  so   
                           $\omega(n-1) $ $\le$ $\gamma$($\omega(n-1) - \omega(n)$)

                           then ($\gamma - 1$)$\omega(n-1)$ $\ge$ $\gamma \omega(n)$ , so $\omega(n)$  $\le$ $\delta$ $\omega(n-1)$,  for $\delta$ = $\frac{\gamma-1}{\gamma}$ $<1$. We repeat this process and we get $\omega(n)$  $\le$ $\delta$ $\omega(n-1)$$\le$ $\delta^{2}$ $\omega(n-2)$ $\le$ $\delta^{3}$ $\omega(n-3)$......
                           
                           finally,$\omega(n)$  $\le$ $\delta^{n}$ $\omega(0)$  for $\omega(0)$ = sup$_{x\in M}$f - inf$_{x\in M}$f.
                           
                           While because  sup$_{x\in M}$ F(x,t) gets smaller when t gets larger, and  inf$_{x\in M}$ F(x,t) gets larger when t gets larger, so $\omega(t)$ decreases when t gets laeger. Therefore, let a = -log($\delta$) , we can choose a constant C$_{4}$ independent of t such that $\omega(t)$ $\le$ C$_{4}$e$^{-at}$.
                           
                           Let $\varphi$(x,t) = $\frac{\partial u}{\partial t}$ - $\frac{1}{Vol(M)}$ $\int_{M}$ $\frac{\partial u}{\partial t}$ $\mathit{d}\tilde{V}$,  Then in order to show when t goes into the infinity, we should analyze the behaviour of $\varphi$(x,t), if we can prove $\varphi$(x,t) goes into 0 then that means $\frac{\partial u}{\partial t}$ truly converges to some function.
                           
                           We use E = $\frac{1}{2}$ $\int_{M} \varphi^{2}$ $\mathit{d}\tilde{V}$  and we want to estimate E in terms of t to figure out how $\varphi$ changes.
                           While 
                           
                            $\mathit{d}$$\tilde{V}$ 
                            
                            =det( $\tilde{\mathit{g_{i\bar{j}}}}$)  $\wedge^{n}_{i = 1}$ ( $\frac{\sqrt{-1}}{2}$ $\mathit{d}z^{i}$ $\wedge$  $\mathit{d}\bar{z}^{j}$) 
                            
                           =det($\mathit{g_{i\bar{j}}}$ +  $\frac{\partial^{2}u}{\partial\mathit{z^{i}}\partial\mathit{\bar{z}^{j}}}$)  $\wedge^{n}_{i = 1}$ ( $\frac{\sqrt{-1}}{2}$ $\mathit{d}z^{i}$ $\wedge$  $\mathit{d}\bar{z}^{j}$)
                            , then 
                            
                            $\frac{\partial}{\partial t}(\mathit{d}\tilde{V})$ 
                            
                            = ($\frac{\partial}{\partial t}$( det($\mathit{g_{i\bar{j}}}$ +  $\frac{\partial^{2}u}{\partial\mathit{z^{i}}\partial\mathit{\bar{z}^{j}}}$))) $\wedge^{n}_{i = 1}$ ( $\frac{\sqrt{-1}}{2}$ $\mathit{d}z^{i}$ $\wedge$  $\mathit{d}\bar{z}^{j}$)
                            
                            =$\tilde{\mathit{g}}^{i\bar{j}}$ $\frac{\partial \tilde{\mathit{g_{i\bar{j}}}}}{\partial t}$ det($\tilde{\mathit{g}}^{i\bar{j}}$) $\wedge^{n}_{i = 1}$ ( $\frac{\sqrt{-1}}{2}$ $\mathit{d}z^{i}$ $\wedge$  $\mathit{d}\bar{z}^{j}$)  
                            
                            =$\frac{\tilde{\mathit{g}}^{i\bar{j}} \frac{\partial \tilde{\mathit{g_{i\bar{j}}}}}{\partial t} det(\tilde{\mathit{g}}^{i\bar{j}})}{det(\tilde{\mathit{g}}^{i\bar{j}}) } det(\tilde{\mathit{g}}^{i\bar{j}}$)  = ($\frac{\partial}{\partial t}$ log det($\mathit{g_{i\bar{j}}}$ +  $\frac{\partial^{2}u}{\partial\mathit{z^{i}}\partial\mathit{\bar{z}^{j}}}$))  $\mathit{d}$$\tilde{V}$
                            
                            =$\frac{\partial}{\partial t}(\frac{\partial u}{\partial t})$ $\mathit{d}$$\tilde{V}$ = $\tilde{\Delta} (\frac{\partial u}{\partial t})$ $\mathit{d}$$\tilde{V}$.
                            
                            Then we calculate 
                            
                            $\frac{\partial\varphi}{\partial t}$(x,t) = $\frac{\partial^{2} u}{\partial t^{2}}(x,t)$ - 
                            $\frac{1}{Vol(M)}$ $\int_{M}$ $\frac{\partial^{2} u}{\partial t^{2}}$ $\mathit{d}\tilde{V}$ - $\frac{1}{Vol(M)}$ $\int_{M}$ $\frac{\partial u}{\partial t}$ $\tilde{\Delta} (\frac{\partial u}{\partial t})$ $\mathit{d}$$\tilde{V}$,
                            
                            while $\frac{1}{Vol(M)}$ $\int_{M}$ $\frac{\partial^{2} u}{\partial t^{2}}$ $\mathit{d}\tilde{V}$ 
                            
                            = 
                            $\frac{1}{Vol(M)}$ $\int_{M}$  $\tilde{\Delta}$ $(\frac{\partial u}{\partial t})$ $\mathit{d}$$\tilde{V}$ 
                            
                            =
                             $\frac{1}{Vol(M)}$ $\int_{M}$  $(\frac{\partial }{\partial t}$ $\mathit{d}$$\tilde{V}$)
                             
                              = 
                             $\frac{\partial }{\partial t}$ ($\frac{1}{Vol(M)}$ $\int_{M}$ $\mathit{d}$$\tilde{V}$) -  $\frac{1}{Vol(M)}$ $\int_{M}$  $(\frac{\partial }{\partial t}$ 1) $\mathit{d}$$\tilde{V}$ = 0-0 = 0,
                             
                             so $\frac{\partial\varphi}{\partial t}$(x,t)   
                             
                            = $\frac{\partial^{2} u}{\partial t^{2}}(x,t)$ - 
                                $\frac{1}{Vol(M)}$ $\int_{M}$ $\frac{\partial u}{\partial t}$ $\tilde{\Delta} (\frac{\partial u}{\partial t})$ $\mathit{d}$$\tilde{V}$ 
                                
                                =
                                 $ \tilde{\Delta}(\frac{\partial u}{\partial t})$- 
                                $\frac{1}{Vol(M)}$ $\int_{M}$ $\frac{\partial u}{\partial t}$ $\tilde{\Delta} (\frac{\partial u}{\partial t})$ $\mathit{d}$$\tilde{V}$.
                                
                                Then we can calculate $\frac{\mathit{d}E}{\mathit{d}t}$  = $\int_{M} \varphi \frac{\partial \varphi}{\partial t}$$\mathit{d}$$\tilde{V}$ 
                                + $\frac{1}{2}$ $\int_{M} \varphi^{2} \tilde{\Delta}(\frac{\partial u}{\partial t})$$\mathit{d}$$\tilde{V}$

                                =$\int_{M}$($\frac{\partial u}{\partial t}$ - $\frac{1}{Vol(M)}$ $\int_{M}$ $\frac{\partial u}{\partial t}$ $\mathit{d}\tilde{V}$)  ( $ \tilde{\Delta}(\frac{\partial u}{\partial t})$- 
                                $\frac{1}{Vol(M)}$ $\int_{M}$ $\frac{\partial u}{\partial t}$ $\tilde{\Delta} (\frac{\partial u}{\partial t})$ $\mathit{d}$$\tilde{V}$)$\mathit{d}$$\tilde{V}$ 
                                
                                + $\frac{1}{2}$ $\int_{M} \varphi^{2} \tilde{\Delta}(\frac{\partial u}{\partial t})$$\mathit{d}$$\tilde{V}$

                                = $\int_{M} $ $\frac{\partial u}{\partial t}$ $ \tilde{\Delta}(\frac{\partial u}{\partial t})$$\mathit{d}$$\tilde{V}$ 
                                
                                - $\int_{M}$ $\frac{1}{Vol(M)}$ $\int_{M}$ $\frac{\partial u}{\partial t}$$\mathit{d}\tilde{V}$($ \tilde{\Delta}(\frac{\partial u}{\partial t})$)$\mathit{d}\tilde{V}$
                                
                                - $\int_{M}$ ($\frac{\partial u}{\partial t}$) ($\frac{1}{Vol(M)}$ $\int_{M}$ $\frac{\partial u}{\partial t}$ $\tilde{\Delta} (\frac{\partial u}{\partial t})$ $\mathit{d}$$\tilde{V}$)  $\mathit{d}\tilde{V}$
                                
                                +$\int_{M}$ ($\frac{1}{Vol(M)}$ $\int_{M}$ $\frac{\partial u}{\partial t}$ $\mathit{d}\tilde{V}$)($\frac{1}{Vol(M)}$ $\int_{M}$ $\frac{\partial u}{\partial t}$ $\tilde{\Delta} (\frac{\partial u}{\partial t})$ $\mathit{d}$$\tilde{V}$) $\mathit{d}\tilde{V}$
                                
                                + $\frac{1}{2}$ $\int_{M} \varphi^{2} \tilde{\Delta}(\frac{\partial u}{\partial t})$$\mathit{d}$$\tilde{V}$ ,
                                
                                while because 
                                $\int_{M}$ ($\frac{1}{Vol(M)}$ $\int_{M}$ $\frac{\partial u}{\partial t}$ $\mathit{d}\tilde{V}$)($\frac{1}{Vol(M)}$ $\int_{M}$ $\frac{\partial u}{\partial t}$ $\tilde{\Delta} (\frac{\partial u}{\partial t})$ $\mathit{d}$$\tilde{V}$) $\mathit{d}\tilde{V}$

                                - $\int_{M}$ $\frac{1}{Vol(M)}$ $\int_{M}$ $\frac{\partial u}{\partial t}$$\mathit{d}\tilde{V}$($ \tilde{\Delta}(\frac{\partial u}{\partial t})$)$\mathit{d}\tilde{V}$
                               
                               - $\int_{M}$ ($\frac{\partial u}{\partial t}$) ($\frac{1}{Vol(M)}$ $\int_{M}$ $\frac{\partial u}{\partial t}$ $\tilde{\Delta} (\frac{\partial u}{\partial t})$ $\mathit{d}$$\tilde{V}$)  $\mathit{d}\tilde{V}$
                               
                               =$\frac{1}{Vol(M)}$ $\int_{M}$ $\frac{\partial u}{\partial t}$ $\mathit{d}\tilde{V}$($\int_{M}$ ($\frac{\partial u}{\partial t}$-1) $\tilde{\Delta} (\frac{\partial u}{\partial t})$ $\mathit{d}$$\tilde{V}$)

                                - $\frac{1}{Vol(M)}$ $\int_{M}$ $\frac{\partial u}{\partial t}$ $\mathit{d}\tilde{V}$ ($\int_{M}$ $\frac{\partial u}{\partial t}$ $\tilde{\Delta} (\frac{\partial u}{\partial t})$ $\mathit{d}$$\tilde{V}$)
                                
                                =$\frac{1}{Vol(M)}$ $\int_{M}$ $\frac{\partial u}{\partial t}$ $\mathit{d}\tilde{V}$($\int_{M}$((-1) $\int_{M}$  $\tilde{\Delta} (\frac{\partial u}{\partial t})$ $\mathit{d}$$\tilde{V}$)
                                = 0
                                
                                 just the same as the arguement before,
                                
                                so $\frac{\mathit{d}E}{\mathit{d}t}$ = $\int_{M} $ $\frac{\partial u}{\partial t}$ $ \tilde{\Delta}(\frac{\partial u}{\partial t})$$\mathit{d}$$\tilde{V}$ 
                                + $\frac{1}{2}$ $\int_{M} \varphi^{2} \tilde{\Delta}(\frac{\partial u}{\partial t})$ $\mathit{d}$$\tilde{V}$, there we use integrating by parts  and M is compact so no boundry  exists, due to $\tilde{\Delta}$ = $\mid \tilde{\nabla} \mid^{2}$, so

                               =- $\int_{M} $ $\tilde{\nabla}$  ($\frac{\partial u}{\partial t}$) $ \tilde{\nabla}(\frac{\partial u}{\partial t})$$\mathit{d}$$\tilde{V}$ 
                                + $\frac{1}{2}$ $\int_{M} $ $ \mathit{d}$$\varphi^{2} \tilde{\nabla}(\frac{\partial u}{\partial t})$ $\mathit{d}$$\tilde{V}$
                                
                                =-$\int_{M} $ $\mid\tilde{\nabla}$   $ \frac{\partial u}{\partial t}\mid^{2}$$\mathit{d}$$\tilde{V}$ 
                                + $\frac{1}{2}$ $\int_{M} $ $ \varphi\mathit{d}$$\varphi \tilde{\nabla}(\frac{\partial u}{\partial t})$ $\mathit{d}$$\tilde{V}$

                                =-$\int_{M} $ $\mid\tilde{\nabla}$   $ \frac{\partial u}{\partial t}\mid^{2}$$\mathit{d}\tilde{V}$
                                +  $\int_{M} $ $ \varphi$$\tilde{\nabla}(\frac{\partial u}{\partial t}) \tilde{\nabla}(\frac{\partial u}{\partial t})$ $\mathit{d}$$\tilde{V}$

                                = $\int_{M} (-1+\varphi) \mid \tilde{\nabla}\frac{\partial u}{\partial t} \mid^{2}$ $\mathit{d}$$\tilde{V}$
                                
                                for $\mid \tilde{\nabla}() \mid^{2}$ = $\tilde{\mathit{g}}^{i\bar{j}}$()$_{i}$ ()$_{j}$ is the square of the gradient.
                                
                                Because sup$_{x\in M}$ $\varphi$(x,t) is the difference between $\varphi$ and its average, while $\omega$ is the largest difference, ans we can choose t large enough such that $\omega$ less than any value, so $\exists t$ such that
                                
                                 sup$_{x\in M}$ $\varphi$(x,t) $<$ $\omega$(t) $<\frac{1}{2}$, then

                                $\frac{\mathit{d}E}{\mathit{d}t}$ 
                                
                                = 
                                $\int_{M} (-1+\varphi) \mid \tilde{\nabla}\frac{\partial u}{\partial t} \mid^{2}$ $\mathit{d}$$\tilde{V}$

                                $\le$ -$\frac{1}{2}$ $\int_{M} \mid \tilde{\nabla}\frac{\partial u}{\partial t} \mid^{2}$ $\mathit{d}$$\tilde{V}$
                                
                                 = -$\frac{1}{2}$ $\int_{M} \mid \tilde{\nabla}\varphi \mid^{2}$ $\mathit{d}$$\tilde{V}$ 
                                
                                And by the definition of $\varphi$,  $\int_{M} \varphi$ $\mathit{d}$$\tilde{V}$ =0,  by the Poincare inequality
                                
                                 $\mid\mid \varphi \mid\mid_{p}$ $\le$ h$\mid\mid D \varphi \mid\mid_{p}$   then $\int_{M} \varphi^{2} \mathit{d}\tilde{V}$ $\le$ h$^{2} \int_{M} \mid D\varphi \mid^{2} \mathit{d}\tilde{V}$ = h$^{2} \int_{M} \mid \tilde{\nabla}\varphi \mid^{2} \mathit{d}\tilde{V}$ for h can be considered as the diameter of the domain,
                                we get 
                                
                                 $\int_{M} \mid \tilde{\nabla}\varphi \mid^{2}$ $\mathit{d}$$\tilde{V}$ $\ge$ $\lambda_{1}(t)$ 
                                $\int_{M} \varphi^{2} $ $\mathit{d}$$\tilde{V}$,  where $\lambda_{1}(t)$ is the first eigenvalue of $\tilde{\Delta}$ at time t.Then there exists a constant C$_{5}$ such that for any t, 
                                
                                $\lambda_{1}(t)$ $>C_{5}$, so  $\frac{\mathit{d}E}{\mathit{d}t}$ = -$\frac{1}{2}$ $\int_{M} \mid \tilde{\nabla}\varphi \mid^{2}$ $\mathit{d}$$\tilde{V}$  $\le$ C$_{5}$E.
                                
                                 That's really similar to an ordinary differential equation, and the exponential function is monotonic, so we can solve the equation and get the inequality :
                                 
                                 E$\le$C$_{6}$ e$^{-C_{5}t}$
                                 
                                 while $\mathit{d}\tilde{V}$ is uniformly equivalent to $\mathit{d}V$, so there exists constant C$^{\prime}_{6}$ such that 
                                
                                $\int_{M} \varphi^{2} \mathit{d}V$ $\le$ C$^{\prime}_{6}$ e$^{-C_{5}t}$,  then we can finally prove the uniform convergence theorem.
                                

                                $\mathbf{Theorem \; of \; uniform\; convergence}$
                                
                                Using the notation in the Proposition of long time existence, as t goes into the infinity, v(x,t) converges to the function v$_{\infty}$ in C$_{\infty}$ topology, therefore, as t goes into the infinity,   
                              $\frac{\partial u}{\partial t}$ converges to a constant in C$_{\infty}$ topology.
                                
                               $ \mathit{Proof.}$
                               
                               First we prove v(x,s) is a Cauchy sequence in L$_{1}$ norm, as t goes into the infinity. For any 0$<s<s^{\prime}$:
                               $\int_{M} \mid v(x,s) - v(x,s^{\prime}) \mid$ $\mathit{d}V$

                               $\le$  $\int_{M} \mid \int_{s}^{s^{\prime}}  \frac{\partial v}{\partial t}(x,t)  \mid \mathit{d}t \mathit{d}V$

                               $\le$  $\int_{M} \int_{s}^{s^{\prime}} \mid \frac{\partial v}{\partial t}(x,t) \mid \mathit{d}t \mathit{d}V$

                               =  $\int_{s}^{s^{\prime}} \int_{M} \mid \frac{\partial v}{\partial t}(x,t) \mid \mathit{d}V \mathit{d}t$

                               =  $\int_{s}^{s^{\prime}} \int_{M} \mid \frac{\partial u}{\partial t} - \frac{1}{Vol(M)} \int_{M} \frac{\partial u}{\partial t} \mathit{d}V \mid \mathit{d}V \mathit{d}t$

                               =$\int_{s}^{s^{\prime}} \int_{M} \mid \frac{\partial u}{\partial t}$ - $\frac{1}{Vol(M)}$ $\int_{M}$ $\frac{\partial u}{\partial t}$$\mathit{d}\tilde{V}$ + $\frac{1}{Vol(M)}$ $\int_{M}$ $\frac{\partial u}{\partial t}$$\mathit{d}\tilde{V}$ 
                               
                               -$\frac{1}{Vol(M)}$ $\int_{M}$ $\frac{\partial u}{\partial t}$$\mathit{d}V$ $\mid \mathit{d}V \mathit{d}t$
                               
                               $\le$ $\int_{s}^{s^{\prime}} \int_{M} \mid \varphi \mid \mathit{d}V \mathit{d}t$
                               + $\int_{s}^{s^{\prime}} \int_{M} \frac{1}{Vol(M)} \mid \int_{M} \frac{\partial u}{\partial t}\mathit{d}\tilde{V} - \int_{M} \frac{\partial u}{\partial t}\mathit{d}V \mid  \mathit{d}V \mathit{d}t$
                               
                               $\le$ $\int_{s}^{s^{\prime}} \int_{M} 1  \mid  \varphi \mid \mathit{d}V \mathit{d}t$
                               + $\int_{s}^{s^{\prime}} \int_{M} \frac{1}{Vol(M)} \mid  sup_{x\in M}(\frac{\partial u}{\partial t})\int_{M}\mathit{d}\tilde{V} $
                               
                               -  $inf_{x\in M}(\frac{\partial u}{\partial t}) \int_{M} \mathit{d}V \mid  \mathit{d}V \mathit{d}t$
                               
                               $\le$ Vol(M)$^{\frac{1}{2}}$ $\int_{s}^{\infty}$($\int_{M} \varphi^{2} \mathit{d}V$)$^{\frac{1}{2}}$$\mathit{d}t$ + Vol(M) $\int_{s}^{\infty} \omega(t)\mathit{d}t$
                               by Cauchy-Schwarz inequality, then by the estimate before we get above term

                                $\le$  Vol(M)$^{\frac{1}{2}}$ $\int_{s}^{\infty}$(C$^{\prime}_{6}$ e$^{-C_{5}t}$)$^{\frac{1}{2}}$$\mathit{d}t$ + Vol(M) $\int_{s}^{\infty} (C_{4}e^{-at})\mathit{d}t$

                                  =C$_{7}$ $\int_{s}^{\infty} e^{-C_{5}\frac{t}{2}}$ $\mathit{d}t$ + C$_{8}$$\int_{s}^{\infty} e^{-at}\mathit{d}t$
                                
                                while the above integral converges in s, so if s goes into the infinity, the two integral can be very small,  so this means v(x,s) is a Cauchy sequence in L$_{1}$ norm. So by L$_{1}$ is complete, there exists a function v$_{\infty}^{\prime}$(x) such that v(x,t) converges uniformly to v$_{\infty}^{\prime}$(x), while in the long time existence theorem we know there exists a sequence t$_{k}$ such that v(x,t$_{k}$) converges to the v$_{\infty}$(x) there, so here v$_{\infty}$(x) =  v$_{\infty}^{\prime}$(x) , so  v(x,t) converges to v$_{\infty}$ in L$_{1}$ norm. However, what we need is the convergence in C$^{\infty}$ topology which can be expressed in the equation we study.
                                
                                We prove it by contradiction. If $\exists$ r$>0$, $\epsilon >0$ such that 
                                
                                $\forall$ N, $\exists$ n$>$N,such that $\mid\mid v(x,t_{n})-v_{\infty}(x) \mid\mid_{C^{r}} >\epsilon$, while N can be chosen arbitary, so we can find a sequence t$_{n}$ such that
                                $\mid\mid v(x,t_{n})-v_{\infty}(x) \mid\mid_{C^{r}}  >\epsilon$, but the sequence v(x,t$_{n}$) is bounded so there exists a subsequence v(x,t$_{kn}$) converges to $\tilde{v}_{\infty}$ $\ne$ v$_{\infty}$(x) in C$^{\infty}$ topology, but we know v(x,t$_{n}$) converges in L$_{1}$ norm to v$_{\infty}$(x),so it's a contradiction. So finally v(x,t) converges to v$_{\infty}$(x) in C$^{\infty}$ topology.
                                
                                Then we consider the equation again:
                                
                                 $\frac{\partial u}{\partial t}$ =
                                log det($\mathit{g_{i\bar{j}}}$ + $\frac{\partial^{2}u}{\partial\mathit{z^{i}}\partial\mathit{\bar{z}^{j}}}$ ) - log det($\mathit{g_{i\bar{j}}}$ ) + f------------------(3.1)
                                
                                because $\partial \tilde{\partial}v$ = $\partial \tilde{\partial}u$, so when t goes into the infinity,
                                
                                 log det($\mathit{g_{i\bar{j}}}$ + $\frac{\partial^{2}u}{\partial\mathit{z^{i}}\partial\mathit{\bar{z}^{j}}}$ ) - log det($\mathit{g_{i\bar{j}}}$ ) + f  converges to
                                 
                                  log det($\mathit{g_{i\bar{j}}}$ + $\frac{\partial^{2}v_{\infty}}{\partial\mathit{z^{i}}\partial\mathit{\bar{z}^{j}}}$ ) - log det($\mathit{g_{i\bar{j}}}$ ) + f, which means  $\frac{\partial u}{\partial t}$ converges to  $\frac{\partial u}{\partial t}_{\infty}(x)$ in C$^{\infty}$ topology when t goes into the infinity, while
                                  
                                   $\omega$(t) = sup$_{x\in M}$$\frac{\partial u}{\partial t}$ - inf$_{x\in M}$$\frac{\partial u}{\partial t}$ $\le$ C$_{4}$e$^{-at}$, when t goes into the infinity, the right hand side goes into 0, therefore, we can only gets $\frac{\partial u}{\partial t}$ converges to a constant.

                                  \section{The final Theorem}

                                  According to our idea:    
                                  First prove the solution exists all the time, thenprove when t goes infinitly, 
                                  $\tilde{\mathit{g_{i\bar{j}}}}$
                                  converges to a definite 
                                  $\tilde{\mathit{g_{i\bar{j}}}}$($\infty$)
                                  and hence 
                                  $\frac{\partial\tilde{\mathit{g_{i\bar{j}}}}}{\partial t}$
                                  converges to 0, then we get

                                  -$\tilde{\mathit{R_{i\bar{j}}}}$ + $\mathit{T_{i\bar{j}}}$ = 0, 
                                  then we get $\mathit{T_{i\bar{j}}}$ 
                                  will be the Ricci tensor of
                                  $\tilde{\mathit{g_{i\bar{j}}}}$($\infty$)
                                  so the $\tilde{\mathit{g_{i\bar{j}}}}$($\infty$) is the metric we want. Then we can finish the whole process.

                                  $\mathbf{Main \; Theorem}$
                                  
                                  M be a compact K$\ddot{a}$hler manifold with the K$\ddot{a}$hler metric  $\mathit{g_{i\bar{j}}}\mathit{d}\mathit{z^{i}}\wedge\mathit{d}\mathit{\bar{z}^{j}}$, C$_{1}$(M) is the first Chern class of M, consider a presentation of it 
                                  $\frac{\sqrt{-1}}{2\pi}\mathit{T_{i\bar{j}}}\mathit{d}\mathit{z^{i}}\wedge\mathit{d}\mathit{\bar{z}^{j}}$, while from the initial metric $\mathit{g_{i\bar{j}}}$, we consider an equation with changing $\mathit{g_{i\bar{j}}}$:
                                  
                                      $\frac{\partial\tilde{\mathit{g_{i\bar{j}}}}}{\partial t}$ = -$\tilde{\mathit{R_{i\bar{j}}}}$ + $\mathit{T_{i\bar{j}}}$ ,  $\tilde{\mathit{g_{i\bar{j}}}}$ = $\mathit{g_{i\bar{j}}}$ at t = 0 
                                      
                                      then the equation exists a long time solution and  $\tilde{\mathit{g}}^{i\bar{j}}$  converges uniformly to another K$\ddot{a}$hler metric $\bar{\mathit{g}}_{i\bar{j}}$ which is in the same K$\ddot{a}$hler class of $\mathit{g_{i\bar{j}}}$ such that
                                      0=-$\bar{\mathit{R_{i\bar{j}}}}$ + $\mathit{T_{i\bar{j}}}$  then

                                      $\bar{\mathit{R_{i\bar{j}}}}$ = $\mathit{T_{i\bar{j}}}$,  which means $\mathit{T_{i\bar{j}}}$ is the Ricci tensor of $\bar{\mathit{g}}^{i\bar{j}}$.
                                      
                                      $\mathit{Proof.}$
                                      
                                      While the de Rahm cohomolgy class of the K$\ddot{a}$hler form  $\frac{\sqrt{-1}}{2\pi}\mathit{R_{i\bar{j}}}\mathit{d}\mathit{z^{i}}\wedge\mathit{d}\mathit{\bar{z}^{j}}$
                                      is the first Chern class C$_{1}$(M) of M, where the $\mathit{R_{i\bar{j}}}$ is the Ricci curvature of the K$\ddot{a}$hler metric. Then because 
                                      	$\frac{\sqrt{-1}}{2\pi}\mathit{T_{i\bar{j}}}\mathit{d}\mathit{z^{i}}\wedge\mathit{d}\mathit{\bar{z}^{j}}$ also represents C$_{1}$(M),
                                      	 so 
                                      	 
                                      $\mathit{T_{i\bar{j}}}$ - $\mathit{R_{i\bar{j}}}$ =  $\frac{\partial^{2}f}{\partial\mathit{z^{i}}\partial\mathit{\bar{z}^{j}}}$,  where f is a real-value smooth function on M. By the long time existence Theorem we know the equation

                                     $\frac{\partial u}{\partial t}$ =
                                     log det($\mathit{g_{i\bar{j}}}$ + $\frac{\partial^{2}u}{\partial\mathit{z^{i}}\partial\mathit{\bar{z}^{j}}}$ ) - log det($\mathit{g_{i\bar{j}}}$ ) + f----------------(3.1)

                                     u(x,t) = 0 when t = 0 exists a smooth solution u(x,t) in all time such that

                                     $\tilde{\mathit{g_{i\bar{j}}}}$(t) - $\mathit{g_{i\bar{j}}}$ = $\frac{\partial^{2}u}{\partial\mathit{z^{i}}\partial\mathit{\bar{z}^{j}}}$
                                      
                                      Then from the uniformly convergence theorem we get that as t goes into the infinity, u(x,t) converges uniformly so $\tilde{\mathit{g}}^{i\bar{j}}$ converges in C$^{\infty}$ topology to  $\tilde{\mathit{g}}^{i\bar{j}}$($\infty$), and by
                                      
                                       $\tilde{\mathit{g_{i\bar{j}}}}$(t) = $\mathit{g_{i\bar{j}}}$ + $\frac{\partial^{2}u}{\partial\mathit{z^{i}}\partial\mathit{\bar{z}^{j}}}$  so
                                     $ \frac{\partial \tilde{\mathit{g}}^{i\bar{j}}}{\partial t}$
                                      converges uniformly to 0 for  $\frac{\partial u}{\partial t}$ converges uniformly.
                                      
                                      Then we differentiate the equation :

                                      $\frac{\partial^{2}}{\partial\mathit{z^{i}}\partial\mathit{\bar{z}^{j}}}$($\frac{\partial u}{\partial t}$) =   
                                      $\frac{\partial^{2}}{\partial\mathit{z^{i}}\partial\mathit{\bar{z}^{j}}}$((log det($\mathit{g_{i\bar{j}}}$ + $\frac{\partial^{2}u}{\partial\mathit{z^{i}}\partial\mathit{\bar{z}^{j}}}$ ) - log det($\mathit{g_{i\bar{j}}}$ )) + $\frac{\partial^{2}f}{\partial\mathit{z^{i}}\partial\mathit{\bar{z}^{j}}}$
                                      
                                       while we can calculate the Ricci tensor explicitly in  terms of
                                       
                                        $\mathit{g_{i\bar{j}}}$ and $\tilde{\mathit{g_{i\bar{j}}}}$:
                                      $\mathit{R_{i\bar{j}}}$ = -$\frac{\partial^{2}}{\partial\mathit{z^{i}}\partial\mathit{\bar{z}^{j}}}$ log det($\mathit{g_{i\bar{j}}}$),  and $ \frac{\partial \tilde{\mathit{g}}^{i\bar{j}}}{\partial t}$  = $\partial \tilde{\partial} (\frac{\partial u}{\partial t})$
                                      
                                      so  $\frac{\partial\tilde{\mathit{g_{i\bar{j}}}}}{\partial t}$ = -$\tilde{\mathit{R_{i\bar{j}}}}$ + $\mathit{T_{i\bar{j}}}$ ,  $\tilde{\mathit{g_{i\bar{j}}}}$ = $\mathit{g_{i\bar{j}}}$ at t = 0  
                                      
                                      then let t goes into the infinity we have 0=-$\bar{\mathit{R_{i\bar{j}}}}$ + $\mathit{T_{i\bar{j}}}$ 
                                      for $\frac{\partial \tilde{\mathit{g}}^{i\bar{j}}(\infty)}{\partial t}$ = 0, it's a constant, so finally
                                      $\bar{\mathit{R_{i\bar{j}}}}(\infty)$ = $\mathit{T_{i\bar{j}}}$. That's the metric we want.
                                      
                                      We prove a corollary as an application.
                                      
                                      $\mathbf{Corollary}$
                                      
                                      If C$_{1}$(M) = 0, thenwe can deform the initial K$\ddot{a}$hler metric to a Ricci flat metric in thr direction of negative Ricci tensor.

                                      $\mathit{Proof.}$
                                      
                                      If C$_{1}M$ = 0, then
                                      $\frac{\partial\tilde{\mathit{g_{i\bar{j}}}}}{\partial t}$ = -$\tilde{\mathit{R_{i\bar{j}}}}$ ,  $\tilde{\mathit{g_{i\bar{j}}}}$ = $\mathit{g_{i\bar{j}}}$ at t = 0  ,
                                      then by the theorem $\tilde{\mathit{g}}^{i\bar{j}}$ converges to $\tilde{\mathit{g}}^{i\bar{j}}$($\infty$), and 0 = -$\tilde{\mathit{R_{i\bar{j}}}}$($\infty$) , which is the Ricci tensor of a Ricci flat metric.

                                        \section{ The problem of existence of the K$\ddot{a}$hler-Einstein metric}
                                        
                                        Here we consider M a compact K$\ddot{a}$hler manifold with negative first Chern class C$_{1}$(M), if we want to find a Kahler-Einstein metric on M i.e. R = kg.  We consider the evolution function:
                                        
                                         $\frac{\partial\tilde{\mathit{g_{i\bar{j}}}}}{\partial t}$ = -$\tilde{\mathit{R_{i\bar{j}}}}$ - $\tilde{\mathit{g}}^{i\bar{j}}$,  $\tilde{\mathit{g_{i\bar{j}}}}$ = $\mathit{g_{i\bar{j}}}$ at t = 0 ,
                                         and   $\mathit{g_{i\bar{j}}}$ is positive definite and represents the negative of the first Chern class, here we know that the first Chern class is the de Rahm cohomology class of the Ricci form which is a real (1,1) form represented as: $\rho$(X,Y) = $\frac{1}{2}$ Rc(JX,Y), and the first Chern class is negative if the tensor Rc is negative definite. Here we choose this evolution equation such that when t goes into the infinity, the Ricci tensor is negative definite for the K$\ddot{a}$hler metric is positive definte by the initial condition, that fits our assumption,  then we derive the scalar equation:
                                         
                                         $\tilde{\mathit{g_{i\bar{j}}}}$(t) = $\mathit{g_{i\bar{j}}}$ + $\frac{\partial^{2}u}{\partial\mathit{z^{i}}\partial\mathit{\bar{z}^{j}}}$

                                         and  $\frac{\partial \tilde{\mathit{g_{i\bar{j}}}}(t) }{\partial t}$ = $\frac{\partial^{2}}{\partial\mathit{z^{i}}\partial\mathit{\bar{z}^{j}}}$($\frac{\partial u}{\partial t}$)
                                         
                                         $\frac{\partial^{2}}{\partial\mathit{z^{i}}\partial\mathit{\bar{z}^{j}}}$($\frac{\partial u}{\partial t}$) =   
                                         $\frac{\partial^{2}}{\partial\mathit{z^{i}}\partial\mathit{\bar{z}^{j}}}$log det($\mathit{g_{i\bar{j}}}$ + $\frac{\partial^{2}u}{\partial\mathit{z^{i}}\partial\mathit{\bar{z}^{j}}}$ ) - $\mathit{g_{i\bar{j}}}$- $\frac{\partial^{2}}{\partial\mathit{z^{i}}\partial\mathit{\bar{z}^{j}}}$ u 
                                         
                                         while $\mathit{g_{i\bar{j}}}$ represents the first Chern class of M, so $\mathit{g_{i\bar{j}}}$ = -$\mathit{R_{i\bar{j}}}$
                                         
                                         so the equation changes to:
                                         
                                         $\frac{\partial^{2}}{\partial\mathit{z^{i}}\partial\mathit{\bar{z}^{j}}}$($\frac{\partial u}{\partial t}$) =   
                                         $\frac{\partial^{2}}{\partial\mathit{z^{i}}\partial\mathit{\bar{z}^{j}}}$log det($\mathit{g_{i\bar{j}}}$ + $\frac{\partial^{2}u}{\partial\mathit{z^{i}}\partial\mathit{\bar{z}^{j}}}$ ) - $\frac{\partial^{2}}{\partial\mathit{z^{i}}\partial\mathit{\bar{z}^{j}}}$log det($\mathit{g_{i\bar{j}}}$)  - $\frac{\partial^{2}u}{\partial\mathit{z^{i}}\partial\mathit{\bar{z}^{j}}}$
                                          and similary :
                                          
                                          $\frac{\partial u}{\partial t}$ =   
                                          log det($\mathit{g_{i\bar{j}}}$ + $\frac{\partial^{2}u}{\partial\mathit{z^{i}}\partial\mathit{\bar{z}^{j}}}$ ) - log det($\mathit{g_{i\bar{j}}}$)  - u + f

                                         Differentiating the equation we get   $\frac{\partial}{\partial t}$($\frac{\partial u}{\partial t}$)
                                         = $\tilde{\Delta} (\frac{\partial u}{\partial t})$ - $\frac{\partial u}{\partial t}$
                                         
                                         which means that 
                                         
                                         $\frac{\partial}{\partial t}$($e^{t}\frac{\partial u}{\partial t}$)
                                         = $\tilde{\Delta} (e^{t}\frac{\partial u}{\partial t})$ 
                                         
                                         that's actually a heat equation and  we use the maximum principle again , which means 
                                         
                                         $\mid \frac{\partial u}{\partial t} \mid$ $\le$ $\mid\frac{\partial u}{\partial t}\mid_{t = 0}$  for when t=0 it's actually a boundary of M$\times$[0,$\infty$), then
                                         
                                         $\mid\frac{\partial u}{\partial t}\mid_{t = 0}$ =   
                                         $\mid$log det($\mathit{g_{i\bar{j}}}$ + $\frac{\partial^{2}u(0)}{\partial\mathit{z^{i}}\partial\mathit{\bar{z}^{j}}}$ ) - log det($\mathit{g_{i\bar{j}}}$)  - u(0) + f$\mid$ =  $\mid f \mid$ is bounded in C$^{0}$, 
                                         then there exists a constant C $>$0 such that  $\mid e^{t}\frac{\partial u}{\partial t} \mid$ $\le$ C  and then 
                                         
                                          $\mid \frac{\partial u}{\partial t} \mid$ $\le e^{-t}$ C which means the exponential decay of the $\frac{\partial u}{\partial t}$.
                                         
                                              For s,t$\ge$0, x$\in$M,
                                              
                                         $\mid u(x,s)-u(x,t) \mid$
                                         
                                          = $\mid \int_{t}^{s} \frac{\partial u}{\partial t}(x,m) \mathit{d}m \mid$

                                         $\le$ $ \int_{t}^{s} \mid\frac{\partial u}{\partial t}(x,m) \mid \mathit{d}m $  
                                         
                                         $\le$ $ \int_{t}^{s} Ce^{-t}\mathit{d}m$ 
                                         
                                         = $C(e^{-t}-e^{-s})$
                                         
                                         which means u(x,t) is a Cauchy sequence in the c$^{0}$(M), then u(t) converges uniformly in C$^{0}$(M)  to some continuous function $u_{\infty}$ on M with a difference of exponential decay,  then we get $\mid\mid u(t)-u_{\infty} \mid\mid_{C^{0}(M)}$ $\le$ $Ce^{-t}$, then the convergent function sequence u(t) is uniformly bounded for t$\in$[0,$\infty$), i.e. $\mid\mid\ u(t) mid\mid_{C^{0}(M)}$ is uniformly bounded for t$\in$[0,$\infty$). Then from above arguement, we obtain the estimates below:
                                         
                                         $\mathbf{Lemma \; 6.1}$

                                         (1).$\exists$ uniform constant C such that $\forall $ t$\in$[0,$\infty$),
                                         
                                          $\mid\mid \frac{\partial u}{\partial t} \mid\mid$ $\le e^{-t}$ C,
                                          
                                          (2).$\exists$ a continuous real-valued function $u_{\infty}$ on M such that $\forall $ t$\in$[0,$\infty$),
                                          
                                          $\mid\mid u(t)-u_{\infty} \mid\mid_{C^{0}(M)}$ $\le e^{-t}$ C,
                                          
                                          (3).$\mid\mid\ u(t) \mid\mid_{C^{0}(M)}$ is uniformly bounded for t$\in$[0,$\infty$)

                                         Here log det($\mathit{g_{i\bar{j}}}$ + $\frac{\partial^{2}u}{\partial\mathit{z^{i}}\partial\mathit{\bar{z}^{j}}}$ ) = -$\frac{\partial u}{\partial t}$    
                                          - log det($\mathit{g_{i\bar{j}}}$)  - u + f 
                                          and from estimate above we claim the RHS is uniformly bounded then the $\tilde{\mathit{g}}^{i\bar{j}}$ is uniformly bounded with a upper bound and a lower bound such that:
                                          there exists a uniformly constant C such that  on M$\times$[0,$\infty$)  and
                                          
                                           $\frac{1}{C}$$\mathit{g}_{i\bar{j}}$ $\le$ $\tilde{\mathit{g}}_{i\bar{j}}$  $\le$ C$\mathit{g}_{i\bar{j}}$, then we can actually prove the final theorem of the existence of the K$\ddot{a}$hler-Ricci flow, but the proof of this estimate refers to lots of canonical estimate of the general K$\ddot{a}$hler Ricci flow, hence in the following sections, we will give more estimates to complete the method.

                                         \section{More  details  from  K$\ddot{a}$hler-Ricci Flow}
                                         
                                         In this section we will figure out some details and problems appeared in the formmer sections, mainly use some other estimates to figure the problem. First we prove the maximal existence time for the K$\ddot{a}$hler Ricci flow equation, then similar as before, we will give some important estimates for the normalized Mange-Ampere equation and the corresponding K$\ddot{a}$hler metric, finally we divide the problem into three condition: the first Chern class positive, negative, or equal to zero, and justify the long time convergence.

                                        	We consider the K$\ddot{a}$hler Ricci flow:
                                        	
                                        	$\frac{\partial}{\partial t}\omega$ = -Ric($\omega$), $\omega = \omega_{0}$ when t = 0
                                        	
                                        	where $\omega$ is a family of K$\ddot{a}$hler form on M, we will note $\omega(t)$ as the solution of the equation. Then we take the cohomology class of the both sides:
                                        	
                                        	$\frac{\partial}{\partial t}[\omega]$ = -C$_{1}(M)$, [$\omega] = [\omega_{0}$] when t = 0
                                         
                                           That is an ordinary differential equation because we assume C$_{1}$(M) is known, and the great mathematician S.S.Chern has proved the C$_{1}$(M) is independent of the metric chosen. We solve it and get :
                                         $[\omega](t)$ = $[\omega_{0}]$ - tC$_{1}(M)$, then if the solution exists, we should get 
                                         
                                         $[\omega_{0}]$ - tC$_{1}(M)$= $[\omega](t)>$0 because the metric is positive definite, then we set
                                         
                                          T = sup\{t$>$0$\mid$$[\omega_{0}]$ - tC$_{1}(M)$$>$0\}, we will prove in [0,T), the solution exists. We note that there $\omega$ = $\sqrt{-1}\mathit{g_{i\bar{j}}}\mathit{d}\mathit{z^{i}}\wedge\mathit{d}\mathit{\bar{z}^{j}}$, 
                                          
                                          and $\omega^{n}$=$n!(\sqrt{-1})^{n}detg$$\mathit{d}\mathit{z^{1}}\wedge\mathit{d}\mathit{\bar{z}^{1}}\wedge$......$\wedge\mathit{d}\mathit{z^{n}}\wedge\mathit{d}\mathit{\bar{z}^{n}}$, so 
                                          
                                          -$\sqrt{-1}\partial\bar{\partial}logdetg$ = -$\sqrt{-1}\partial\bar{\partial}logdet\omega^{n}$.
                                         
                                         Let $\eta\in$ $[\omega_{0}]$ - T$^{\prime}$C$_{1}(M)$ be in a K$\ddot{a}$hler class, then let $\hat{\omega}_{t}$=$\frac{1}{T^{\prime}}((T^{\prime}-t)\omega_{0}+t\eta)$, it is still in $[\omega_{0}]$ - T$^{\prime}$C$_{1}(M)$ because it is in the path from one to the other. We assume  	$\frac{\partial}{\partial t}\hat{\omega}_{t}$ = 
                                         $\frac{\sqrt{-1}}{2\pi}\partial\bar{\partial}log\Omega$ in the -C$_{1}$(M), existence of this $\Omega$ is from [9]. Then we find that the solution $\omega$(t) is equivalent to $\omega(t)$=$\hat{\omega}_{t}$+$\frac{\sqrt{-1}}{2\pi}\partial\bar{\partial}\tilde{\varphi}(t)$ exists using the $\partial\tilde{\partial}$  Lemma and $\int_{M}\tilde{\varphi(t)}\omega_{0}^{n}$=0, such that $\tilde{\varphi} $ is smooth on M$\times$[0,T$^{\prime}$) by  the regularity theorem so the equation is  
                                         
                                         $\frac{\sqrt{-1}}{2\pi}\partial\bar{\partial}log\omega^{n}$=$\frac{\partial}{\partial t}$($\hat{\omega}_{t}$+$\frac{\sqrt{-1}}{2\pi}\partial\bar{\partial}\tilde{\varphi}(t)$) = $\frac{\sqrt{-1}}{2\pi}\partial\bar{\partial}log\Omega$ + $\frac{\sqrt{-1}}{2\pi}\partial\bar{\partial}$($\frac{\partial\tilde{\varphi}}{\partial t}$);
                                         then
                                         
                                           $\frac{\sqrt{-1}}{2\pi}\partial\bar{\partial}(log\omega^{n}$ - $log\Omega$ - $\frac{\partial\tilde{\varphi}}{\partial t}$)=0;
                                           
                                           which means 
                                           
                                           c(t)+ $log\frac{\omega^{n}}{\Omega}$  = $\frac{\partial\tilde{\varphi}}{\partial t}$;
                                           
                                           where c(t) is a function just for t. Then  let $\varphi(t)$ = $\tilde{\varphi}(t)$- $\int_{0}^{t}c(s)\mathit{d}s-\tilde{\varphi(0)}$, we get this $\varphi(t)$ solves the equation 
                                      
                                        $\frac{\partial\tilde{\varphi}}{\partial t}$=$log\frac{\omega^{n}}{\Omega}$ for $\varphi(0)$=0,
                                        
                                        conversely, if this $\varphi(t)$ solves the equation in the [0,T$^{\prime}$)
                                        
                                        $\frac{\partial\tilde{\varphi}}{\partial t}$=$log\frac{\omega^{n}}{\Omega}$ for $\varphi(0)$=0,
                                        then let $\omega(t)$ = $\hat{\omega}_{t}$+$\frac{\sqrt{-1}}{2\pi}\partial\bar{\partial}\varphi(t)$,
                                        then
                                        
                                         $\frac{\partial\omega}{\partial t}$ 
                                         
                                         = $\frac{\partial \hat{\omega}_{t}}{\partial t}$ + $\frac{\sqrt{-1}}{2\pi}\partial\bar{\partial}(\frac{\partial}{\partial t}\varphi(t))$
                                         
                                         =$\frac{\sqrt{-1}}{2\pi}\partial\bar{\partial}log\Omega$+
                                        $\frac{\sqrt{-1}}{2\pi}\partial\bar{\partial}$($log\frac{\omega^{n}}{\Omega}$)
                                        
                                         = 
                                      $\frac{\sqrt{-1}}{2\pi}\partial\bar{\partial}$($log\Omega$+$log\frac{\omega^{n}}{\Omega}$)
                                      
                                       = 
                                  $\frac{\sqrt{-1}}{2\pi}\partial\bar{\partial}$logdet$\omega$=-Ric($\omega$),  $\omega(0)=\omega_{0}$, so the existence of $\varphi$ guarantee the existence of the solution of the original equation. Remark if the solution exits, section has proved the uniqueness of the solution. Just the same as before, we should give estimates to this $\varphi$ and first give the short time existence, but we will assume a T$_{max}$ such that there is a solution only in [0,T$_{max}$),  for T$_{max}<$T, then by finding a solution in [T$_{max}$,T) conclude a contradiction then prove the maximum of T.
                                      
                                      $\frac{\partial\varphi}{\partial t}$=$log\frac{(\hat{\omega}_{t}+\frac{\sqrt{-1}}{2\pi}\partial\bar{\partial}\varphi(t))^{n}}{\Omega}$ for $\varphi(0)$=0.
                                      
                                      First we need a lemma.
                                      
                                     $ \mathbf{Lemma \; 7.1}$
                                     
                                     T$>0$,f(x,t) a smooth function on M$\times$[0,T], if f attains its maximum(minimum) at point (x$_{0}$,$t_{0}$), then either $t_{0}$=0 or at (x$_{0}$,$t_{0}$) we have:
                                      	$\frac{\partial f}{\partial t}$$\ge0(\le0)$, $\mathit{d}f=0$, $\sqrt{-1}\partial\tilde{\partial}f\le0(\ge0)$.
                                      	
                                      	$\mathit{Sketch \; of \; Proof.}$ We know that if a smooth function attains its maximum at a point (x$_{0}$,$t_{0}$), then it has zero first derivative and nonpositive Hessian at this point from high dimensional Taylor expension, if $t_{0}>0$, then f is nondecreasing at  $t_{0}$, so $\frac{\partial f}{\partial t}\ge0$, while if $t_{0}=0$, because it is the maximum, so we can only attains f is nonpositive from $t_{0}$.
                                     
                                      Back to the estimate.We should estimate the zero order of the solution $\varphi(T)$. The arguement is due to Song and Weinkove[9].
                                      
                                      $\mathbf{Proposition \; 7.2}$
                                      
                                      $\exists$ C constant such that $\forall$t$\in$[0,T$_{max}$), $\mid\mid\varphi(t)\mid\mid_{C^{0}(M)}\le C$.
                                      
                                      $\mathit{Proof.}$
                                      
                                      Consider $\theta(t)$ = $\varphi(t)-At$, where A is to be determined, because $\varphi(0)=0$, so if we can choose an A to prove $\theta$ attains its maximum at t = 0, then we can give a uniform estimate of $\varphi(t)$. 
                                      
                                      The $\theta$ satisfy the equation $\frac{\partial\theta}{\partial t}$=$log\frac{(\hat{\omega}_{t}+\frac{\sqrt{-1}}{2\pi}\partial\bar{\partial}\theta(t))^{n}}{\Omega}$-A
                                      in [0,t$^{\prime}$]
                                      
                                       for $<0t^{\prime}<T_{max}$. We assume $\theta$ attains its maximum at (x$_{0},t_{0}$) in M$\times$[0,t$^{\prime}$] compact, then we use the Lemma 7.1, if $t^{0}>0$ then
                                      $\frac{\partial \theta}{\partial t}$$\ge0$ and $\sqrt{-1}\partial\tilde{\partial}\theta\le0$, so

                                      0$\le$$\frac{\partial\theta}{\partial t}$
                                      
                                      =$log\frac{(\hat{\omega}_{t}+\frac{\sqrt{-1}}{2\pi}\partial\bar{\partial}\theta(t))^{n}}{\Omega}$-A

                                      $\le$ $log\frac{(\hat{\omega}_{t})^{n}}{\Omega}$-A, if we chosen A large enough, then this inequality can not hold, then we choose a appropriate A such that t$_{0}=0$ is the only choice. Then $\theta(t)\le\theta(0)=0$, which means $\varphi(t)\le At^{\prime}\le AT_{max}$, hence the upper bound.
                                      For the lower bound, we choose $\theta(t)=\varphi(t)+At$, then we consider the minimum of $\theta(t)$, use the Lemma 7.1 again, similarly if $t_{0}>0$ then $\frac{\partial \theta}{\partial t}$$\le0$ and $\sqrt{-1}\partial\tilde{\partial}\theta\ge0$, then

                                      0$\ge$$\frac{\partial\theta}{\partial t}$
                                      
                                      =$log\frac{(\hat{\omega}_{t}+\frac{\sqrt{-1}}{2\pi}\partial\bar{\partial}\theta(t))^{n}}{\Omega}$+A

                                      $\ge$ $log\frac{(\hat{\omega}_{t})^{n}}{\Omega}$+A, 
                                      
                                      then we can chooose A large enough such that the inequality does not hold, hence again we get $\theta(t)\ge\theta(0)=0$, then $\varphi(t)\ge-AT_{max}$, that is the lower bound. We note that for a real (1,1) form $\alpha$=$\frac{\sqrt{-1}}{2\pi}\partial\bar{\partial}\alpha_{i\bar{j}}$$\mathit{d}\mathit{z^{i}}\wedge\mathit{d}\mathit{\bar{z}^{j}}$, the trace with respect to $\omega$ is defined as
                                      
                                       tr$_{\omega}\alpha$= $\mathit{g^{i\bar{j}}}\alpha_{i\bar{j}}$=$\sum$$\frac{\alpha_{i\bar{j}}}{\omega_{i\bar{j}}}$=$\frac{n\omega^{n-1}\wedge\alpha}{\omega^{n}}$, and $\Delta f$=$tr_{\omega}(\frac{\sqrt{-1}}{2\pi}\partial\bar{\partial}f)$, actually this notation we use in the later sections have the same property as the normal trace and laplace.
                                      
                                      Next we want to give an estimate of $\frac{\partial \varphi}{\partial t}$, the arguement refers to Tian-Zhang.
                                      
                                      $\mathbf{Proposition \; 7.3}$
                                      
                                      $\exists$ C positive constant such that on M$\times$[0,T$_{max}$),
                                      
                                       $\frac{1}{C}\Omega\le \omega^{n}(t)\le C\Omega$, equivalently, $\mid\mid \frac{\partial \varphi}{\partial t}\mid\mid_{C^{0}}$ is uniformly bounded.
                                       
                                       $\mathit{Proof.}$
                                       
                                       $\frac{\partial}{\partial t}$$ log\frac{\omega^{n}(t)}{\omega^{n}(0)}$

                                       =$\frac{\partial}{\partial t}$(log$\omega^{n}(t)$-log$\omega^{n}(0)$)

                                       =$\frac{\partial}{\partial t}$(log det$\mathit{g_{i\bar{j}}}$)

                                       =$\mathit{g^{i\bar{j}}}$ 	$\frac{\partial}{\partial t}$ $\mathit{g_{i\bar{j}}}$, 
                                       
                                       but we take trace of the both sides of the equation 	
                                       $\frac{\partial}{\partial t}\omega$ = -Ric($\omega$);
                                       we get $\mathit{g^{j\bar{i}}}$ 	$\frac{\partial}{\partial t}$   $\mathit{g_{i\bar{j}}}$=-R,  where R is the scalar curvature, then  
                                       $\frac{\partial}{\partial t}$$ log\frac{\omega^{n}(t)}{\omega^{n}(0)}$ 
                                       =-R.
                                     
                                       While R =tr(Ric), so 
                                       
                                       R = - $\mathit{g^{\bar{j}i}}$ $\partial_{i}\bar{\partial_{j}}log det g$, 
                                       
                                       then 
                                       $\frac{\partial}{\partial t}$ R 
                                       
                                       = 
                                       - $\mathit{g^{\bar{j}i}}$ $\partial_{i}\bar{\partial_{j}}$($\mathit{g^{l\bar{k}}}$ 	$\frac{\partial}{\partial t}$   $\mathit{g_{k\bar{l}}}$)-	$\frac{\partial}{\partial t}$$\mathit{g^{i\bar{j}}}$ $\partial_{i}\bar{\partial_{j}}log det g$

                                      =- $\mathit{g^{\bar{j}i}}$ $\partial_{i}\bar{\partial_{j}}$(-R)- $\mathit{g^{\bar{l}i}}$ $\mathit{g^{\bar{j}k}}$ $\mathit{R_{k\bar{l}}}$ $\mathit{R_{i\bar{j}}}$, 
                                     
                                      here we calculate $\frac{\partial}{\partial t}$$\mathit{g^{i\bar{j}}}$ by derivative  
                                      $\mathit{g^{i\bar{j}}}$$\mathit{g_{i\bar{j}}}$ = $\delta^{i}_{j}$.
                                      
                                      Then because R =tr(Ric), hence n$\mid Ric(\omega) \mid^{2}\ge R^{2}$, so  (	$\frac{\partial}{\partial t}$ - $\Delta$)R $\ge$ $\frac{1}{n}$R$^{2}$, 
                                      
                                      hence 
                                       ($\frac{\partial}{\partial t}$ - $\Delta$)R $\ge$0, by the maximum principle $\mid R\mid \le C$ .
                                       Integrating  $\frac{\partial}{\partial t}$$ log\frac{\omega^{n}(t)}{\omega^{n}(0)}$=-R , then we get
                                       $\omega^{n}(t)\le e^{Ct}\omega^{n}(0)$. So we finish the upper bound. As for the lower bound, take the derivative of the equation

                                      $\frac{\partial\varphi}{\partial t}$=$log\frac{(\hat{\omega}_{t}+\frac{\sqrt{-1}}{2\pi}\partial\bar{\partial}\varphi(t))^{n}}{\Omega}$ for $\varphi(0)$=0.
                                      
                                      Then we get 
                                      
                                       $\frac{\partial^{2} \varphi}{\partial t^{2}}$ 
                                       
                                       = $\frac{\Omega}{\omega^{n}}$n$\frac{\omega^{n-1}}{\Omega}$ ($\frac{\partial \hat{\omega}_{t}}{\partial t}$+$\frac{\sqrt{-1}}{2\pi}\partial\bar{\partial}$ $\frac{\partial \varphi}{\partial t}$)

                                      =$\Delta \frac{\partial \varphi}{\partial t} $+tr$_{\omega}$$\frac{\partial \hat{\omega}_{t}}{\partial t}$,
                                       and due to  $\Delta \varphi$
                                       
                                       =$tr_{\omega}(\frac{\sqrt{-1}}{2\pi}\partial\bar{\partial}\varphi)$
                                       
                                       =tr$_{\omega}$($\omega-\hat{\omega}_{t}$)
                                       
                                       =n-tr$_{\omega}$$\hat{\omega}_{t}$,
                                       
                                       we let Q = (T$^{\prime}-t$)	$\frac{\partial \varphi}{\partial t}$+$\varphi$+nt, then we calculate 
                                       
                                       $\frac{\partial}{\partial t}$Q = -$\frac{\partial \varphi}{\partial t}$+(T$^{\prime}$-t)$\frac{\partial^{2} \varphi}{\partial t^{2}}$+$\frac{\partial \varphi}{\partial t}$+n,
                                      $\Delta $Q=-(T$^{\prime}$-t)$\Delta $$\frac{\partial \varphi}{\partial t}$+$\Delta \varphi$,
                                      
                                      hence ($\frac{\partial}{\partial t}$-$\Delta$)Q 
                                      
                                      = (T$^{\prime}$-t)tr$_{\omega}$$\frac{\partial \hat{\omega}_{t}}{\partial t}$+n-$\Delta \varphi$
                                      
                                      =tr$_{\omega}$((T$^{\prime}$-t)$\frac{\partial \hat{\omega}_{t}}{\partial t}$+  $\hat{\omega}_{t}$)
                                      
                                      =tr$_{\omega}$((T$^{\prime}$-t) $\frac{1}{T^{\prime}}(\eta-\omega_{0})$+$\frac{1}{T^{\prime}}$((T$^{\prime}$-t)$\omega_{0}+t\eta$))

                                      =tr$_{\omega}$($\frac{1}{T^{\prime}}$t$\eta$)=tr$_{\omega}$$\hat{\omega}_{T^{\prime}}>0$  because it is positive definite.
                                      Then hence for this elliptic equation we use the maximum principle, Q is larger then its infimum of the boundry value at t=0 in M$\times$[0,T$_{max}$), i.e.
                                      
                                      (T$^{\prime}-t$)	$\frac{\partial \varphi}{\partial t}$+$\varphi$+nt
                                      
                                      = Q$\ge$ $\mid T^{\prime} inf_{M} log \frac{\omega_{0}}{\Omega} \mid$, 
                                      
                                      then  $\mid\mid\frac{\partial \varphi}{\partial t}\mid\mid \ge$ $\frac{1}{\mid T^{\prime}-T_{max}\mid}$($\mid$$\mid T^{\prime} inf_{M} log \frac{\omega_{0}}{\Omega} \mid$$\mid$+$\mid\mid \varphi\mid\mid$+nT$_{max}$), by the zero order estimate of $\varphi$ we conclude the lower estimate of $\frac{\partial \varphi}{\partial t}$, note these norms are uniformly C$_{0}$(M) norm.

                                        $\mathbf{Proposition\;7.4}$
                                        
                                        There exists a uniform positive constant C such that on M$\times$[0,T$_{max}$),
                                        
                                        $\frac{1}{C}\omega_{0}$ $\le$  $\omega$ $\le$ $C\omega_{0}$.
                                        
                                        $\mathit{Proof.}$
                                        First we estimate the lower bound. Due to the existence of the normal coordinates, we can always find a system of coordinates for the Hermitian metric g such that: (g$_{0}$)$_{i\bar{j}}$ = $\delta^{i}_{j}$, and $\mathit{g_{i\bar{j}}}$=$\lambda_{i}\delta^{i\bar{j}}$, for positive eigenvalue $\lambda_{i}$.
                                        Then 
                                        
                                        tr$_{\omega}\omega_{0}$=$\sum_{i}^{n}$ $\frac{1}{\lambda_{i}}$
                                        
                                         $\le$ $\frac{1}{(n-1)!}$ $ \frac{(\sum_{i}^{n} \lambda_{i})^{n-1}}{\lambda_{1}...\lambda_{n}}$
                                         
                                         =$\frac{1}{(n-1)!}(tr_{\omega_{0}}\omega)^{n-1}$ $ \frac{\omega_{0}^{n}}{\omega^{n}}\le$C, because from Proposition 7.3 $\omega^{n}$ is uniformly bounded.
                                        
                                        As for the upper bound, we prove  there exists a uniform constant C such that tr$_{\omega_{0}}\omega\le$C on M $\times$[0,T$_{max}$). We use the traditional method: find a useful quantity using coeffecients to be determined,then we prove it satisfies an elliptice equation and use the maximum principle, that method has benn used in many of our proof, but how to find a useful quantity is essential, this can only be attained by experience and attampt. We consider Q = log tr$_{\omega_{0}}\omega$ -A$\varphi$, where A is a positive uniform constant to be determined.
                                      We fix t$^{\prime}\in$(0,T$_{max}$), and assume Q attains its maximum at point (x$_{0}$,t$_{0}$), with out loss generality let t$_{0}>0$.
                                      
                                      First we estimate (	$\frac{\partial}{\partial t}$-$\Delta$)tr$_{\omega_{0}}\omega$.
                                      Let$\hat{\omega}$ be a fixed K$\ddot{a}$hler form corresponding to the K$\ddot{a}$hler-Ricci metric $\hat{g}$ on M, and the same as before, $\omega$ the solution of the original K$\ddot{a}$hler-Ricci flow equation, we note R$_{i\bar{j}}^{\bar{q}p}$=$\mathit{g^{\bar{q}p}}$$\mathit{g^{\bar{l}p}}$R$_{i\bar{j}k\bar{l}}$ is the curvature, $\nabla$ is the connection corresponding to g, $\hat{R}$$_{i\bar{j}}^{\bar{q}p}$=$\mathit{\hat{g}^{\bar{q}p}}$$\mathit{\hat{g}^{\bar{l}p}}$$\hat{R}$$_{i\bar{j}k\bar{l}}$ is the curvature, $\hat{\nabla}$ is the connection corresponding to $\hat{\omega}$. Then the same as before, we use the normal coordinates for $\hat{g}$, then

                                       $\Delta$tr$_{\hat{\omega}}\omega$ 
                                       
                                       = $g^{\bar{l}k}$ $\partial_{k}$ $\partial_{\bar{l}}$($\hat{g}^{\bar{j}i}$ $g_{i\bar{j}}$)
                                       
                                       =
                                      $g^{\bar{l}k}$ ($\partial_{k}$ $\partial_{\bar{l}}$$\hat{g}^{\bar{j}i}$) $g_{i\bar{j}}$
                                      + $g^{\bar{l}k}$ $\hat{g}^{\bar{j}i}$ ($\partial_{k}$ $\partial_{\bar{l}}$ $g_{i\bar{j}}$)

                                      =$g^{\bar{l}k}$ R$_{j\bar{i}}^{k\bar{l}}$ $g_{i\bar{j}}$ - $\hat{g}^{\bar{j}i}$R$_{i\bar{j}}$ 
                                      + $\hat{g}^{\bar{j}i}$ g$^{\bar{q}p}$  g$^{\bar{l}k}$ $\partial_{i}$ g$_{p\bar{l}}$ $\partial_{\bar{j}}$ g$_{k\bar{q}}$, 
                                      
                                       while  $\frac{\partial}{\partial t}$ tr$_{\hat{\omega}} \omega$ = - $\hat{g}^{\bar{j}i}$R$_{i\bar{j}}$, 
                                       
                                       we get
                                        (	$\frac{\partial}{\partial t}$-$\Delta$)tr$_{\hat{\omega}} \omega$ = $g^{\bar{l}k}$ R$_{j\bar{i}}^{k\bar{l}}$ $g_{i\bar{j}}$ - $\hat{g}^{\bar{j}i}$ $g^{\bar{q}p}$  $g^{\bar{l}k}$ $\hat{\nabla}_{i}$ $g_{p\bar{l}}$ $\hat{\nabla}_{\bar{j}}$ $g_{k\bar{q}}$.
                                       
                                      Then we use the formular above to calculate (	$\frac{\partial}{\partial t}$-$\Delta$)log tr$_{\hat{\omega}} \omega$, then choose normal coordinates such that g is diagnoal, by Cauchy-Schwarz inequality we get
                                      
                                       $\mid \partial tr_{\hat{\omega}} \omega \mid^{2}_{g}$ $\le$
                                      (tr$_{\hat{\omega}} \omega$) $\sum_{i,j,k}^{}$  $g^{\bar{i}i}$  $g^{\bar{j}j}$ $\partial_{k}$ $g_{i\bar{j}}$ $\partial_{\bar{k}}$ $g_{j\bar{i}}$,

                                      and then we define $\hat{C}$ = -inf$_{x\in M}${$\hat{R}_{i\bar{i}j\bar{j}}$(x)$\mid $ {$\partial_{z^{1}}$,...,$\partial_{z^{n}}$ }}, then calculate
                                      
                                       $g^{\bar{l}k}$ $\hat{R}^{j\bar{i}}_{k\bar{l}}$ $g_{i\bar{j}}$
                                       
                                       =$\sum_{k,i}^{}$ $g^{\bar{k}k}$ $\hat{R}_{k\bar{k}i\bar{i}}$ $g_{i\bar{i}}$ 
                                       
                                       $\ge$ -$\hat{C}$ $\sum_{k}$ $g^{\bar{k}k}$ $\sum_{i}$ $g^{i\bar{i}}$
                                       
                                        = -$\hat{C}$(tr$_{\hat{\omega}} \omega$)(tr$_{\hat{\omega}} \omega$), then use the equality above 
                                      (	$\frac{\partial}{\partial t}$-$\Delta$)log tr$_{\hat{\omega}} \omega$ $\le$ $\hat{C}$ tr$_{\hat{\omega}} \omega$. This proof comes from [9], one can obtains more details from [9].

                                      Then at point ($x_{0}$,t$_{0}$) we use the estimate above.
                                      
                                      0$\le$($\frac{\partial \varphi}{\partial t}$-$\Delta$)Q
                                      
                                      $\le$ C$_{0}$tr$_{\omega}\omega_{0}$-A	$\frac{\partial \varphi}{\partial t}$+A$\Delta \varphi$ 
                                      
                                      = tr$_{\omega}$($C_{0}\omega_{0}-A\hat{\omega}_{t_{0}}$)-Alog$\frac{\omega^{n}}{\Omega}$+An, 
                                      
                                      do not forget

                                       $\Delta \varphi$
                                       
                                       =
                                       $tr_{\omega}(\frac{\sqrt{-1}}{2\pi}\partial\bar{\partial}\varphi)$

                                       =tr$_{\omega}$($\omega-\hat{\omega}_{t}$)
                                       
                                       =n-tr$_{\omega}$$\hat{\omega}_{t}$, for C$_{0}$ is a constant only depends on the lower bound of the bisectional curvature of g$_{0}$.
                                       Then we need choose A large enough such that
                                       
                                        A $\hat{\omega}_{t_{0}}$-($C_{0}$+1)$\omega_{0}$ is K$\ddot{a}$hler on M, i.e. first it should be positive definite. Then 
                                       
                                       tr$_{\omega}$(A $\hat{\omega}_{t_{0}}$-($C_{0}$+1)$\omega_{0}$) $\ge0$, so
                                       
                                        tr$_{\omega}$(-A $\hat{\omega}_{t_{0}}$+$C_{0}$$\omega_{0}$) $\le$ -tr$_{\omega}\omega_{0}$,
                                      
                                       when at point (x$_{0},t_{0}$),

                                       0$\le$ tr$_{\omega}$($C_{0}\omega_{0}-A\hat{\omega}_{t_{0}}$)-Alog$\frac{\omega^{n}}{\Omega}$+An, 
                                       
                                       so 
                                       tr$_{\omega}$$\omega_{0}$+Alog$\frac{\omega^{n}}{\Omega}$ $\le$An, 
                                       
                                       then  tr$_{\omega}$$\omega_{0}$+Alog$\frac{\omega^{n}}{\omega_{0}}$ $\le$C, for a constant C.
                                      Similarly, find a system of coordinates for the Hermitian metric g such that: (g$_{0}$)$_{i\bar{j}}$ = $\delta^{i}_{j}$, and $\mathit{g_{i\bar{j}}}$=$\lambda_{i}\delta^{i\bar{j}}$, for positive eigenvalue $\lambda_{i}$, then it shows 
                                      
                                      $\sum_{i}^{n}$$\frac{1}{\lambda_{i}}$ + $\sum_{i}^{n}$Alog$\lambda_{i}$ $\le$C. Due to $\lambda_{i}$ is positive, so there is a uniform upper bound C for $\frac{1}{\lambda_{i}}$+Alog$\lambda_{i}$ $\le$C, while $\frac{1}{\lambda_{i}}$ is positive, we attains a uniform upper bound for Alog$\lambda_{i}$, hence a uniform upper bound for $\lambda_{i}$, i.e. $\lambda_{i}\le C$, then we attains an upper bound for tr$_{\omega_{0}}\omega$ at (x$_{0}$,t$_{0}$), then by the uniform bounded $\varphi$  on M$\times$[0,t$^{\prime}$) when t$^{\prime}\le T_{max}$and tr$_{\omega_{0}}\omega$ at (x$_{0}$,$t_{0}$), Q is uniformly bounded on  M$\times$[0,t$^{\prime}$) when t$^{\prime}\le T_{max}$,then use the uniformly bounded $\varphi$ on $[0,T_{max}]$ , we know that tr$\omega_{0}\omega$ has a uniformly upper bound. Therefore, we have prove the uniformly bound for 
                                       $\frac{\partial\varphi}{\partial t}$=$log\frac{\omega^{n}}{\Omega}$ .
                                       
                                       $\mathbf{Lemma \;7.5}$
                                       
                                       If the solution $\omega(t)$ of the K$\ddot{a}$hler-Ricci flow equation before on M$\times[0,T)$ satisfies there $\exists$ a constant C$_{0}$ such that :
                                       $\frac{1}{C_{0}}\omega_{0}$ $\le$ $\omega$ $\le$ $C_{0}\omega_{0}$.
                                       
                                       Then for any positive integer m, there exists corresponding uniform constant $C_{m}$ such that $\mid\mid \omega(t) \mid\mid_{C^{m}(g_{0})}$ $\le$ $C_{m}$.
                                       
                                       $\mathit{Proof.}$
                                       
                                       If the solution $\omega(t)$ satisfies there $\exists$ a constant C$_{0}$ such that :
                                       
                                       $\frac{1}{C_{0}}\omega_{0}$ $\le$ $\omega$ $\le$ $C_{0}\omega_{0}$, we prove first there exists constants C, C$^{\prime}$ depending only on C$_{0}$ and $\omega_{0}$ such that:

                                       $\mid \nabla_{g_{0}}g \mid^{2}$ $\le$C, and ($\frac{\partial}{\partial t}$-$\Delta$)$\mid \nabla_{g_{0}}g \mid^{2}$ $\le$ -$\frac{1}{2}$ $\mid R_{i\bar{j}k\bar{l}} \mid^{2}$ + $C^{\prime}$;
                                       
                                       then  there exists constants C, C$^{\prime}$ depending only on C$_{0}$ and $\omega_{0}$ such that:
                                       
                                        $\mid R_{i\bar{j}k\bar{l}} \mid^{2}$ $\le$C and ($\frac{\partial}{\partial t}$-$\Delta$)
                                       $\mid R_{i\bar{j}k\bar{l}} \mid^{2}$ $\le$ -$\mid \nabla R_{i\bar{j}k\bar{l}}  \mid^{2}$-$\mid \bar{\nabla} R_{i\bar{j}k\bar{l}}  \mid^{2}$ +C$^{\prime}$, where $\bar{\nabla}$ is the conjugate of $\nabla$, then using the condition above we learn that $\exists$ uniform constants C$_{m}$ for positive integer m such that $\mid \nabla^{m}_{R} R_{i\bar{j}k\bar{l}} \mid^{2}$ $\le$C$_{m}$, therefore using the conclusions of the above claim, for U an open subsets of M, for any  compact subset K in U, positive integer m, there exists constants C$^{\prime}_{m}$ depending only on $\omega_{0}$, K,U, and C$_{m}$ such that $\mid\mid \omega(t) \mid\mid_{C^{m}(K,g_{0})}$ $\le$C$^{\prime}_{m}$, finally, from the above claim, we conclude the lemma.

                                       Due to the much too long standard proof of this lemma, we give a sketch of proof, the remaining detail one can refer to [9], Theorem 2.13,2.14,2.15.
                                       
                                       Now we prove the existence of the K$\ddot{a}$hler-Ricci flow solution.
                                       
                                       $\mathbf{Theorem \; of \; the\; maximal \; existence \; time}$
                                       
                                       The K$\ddot{a}$hler-Ricci flow equation
                                       
                                       	$\frac{\partial}{\partial t}\omega$ = -Ric($\omega$), $\omega = \omega_{0}$ when t = 0
                                       	
                                       	has a unique solution in the maximal time $t\in[0,T)$, then this solution exists for all time.
                                       
                                       $\mathit{Proof.}$
                                       
                                       By Proposition 7.4 and Lemma 7.5, we conclude the uniform C$^{\infty}$ estimate for $\omega$(t) on [0,T$_{max}$), then as t goes into the T$_{max}$, by Arzela-Ascoli Theorem and take countable diagonal subsequences, we obtain a solution g(T$_{max}$) while g converges to it on [0,T$_{max}$], then we use the same arguement as solving the problem before, our estimates are independent of t, then our estimates still works the estimates above are independent of t ,so if we choose t$_{0}\in[0,T)$ then the solution also exists in [t$_{0}$,t$_{0}+\epsilon$] for $\epsilon$ independent of t, so we can use this  $\tilde{\mathit{g_{i\bar{j}}}}$(T) as the intial condition of the same equation chosing a new initial point and continue to deformation, which is contradict to the definition of T$_{max}$, hence the process can continue to infinity because our estimates always work,  so the solution exists for all time.
                                       
                                       Then we complete the proof of the condition  C$_{1}$(M) $<0$ which remains.

                                       We assume C$_{1}(M)<0$, and [$\omega_{0}$]=-C$_{1}$(M).We consider the normalized K$\ddot{a}$hler-Ricci flow:
                                       
                                       	$\frac{\partial}{\partial t}\omega$ = -Ric($\omega$)-$\omega$, $\omega = \omega_{0}$ when t = 0
                                       	
                                       	We use this normalized form to avoid the K$\ddot{a}$hler class [$\omega$(t)] given by (1+t)[$\omega_{0}$] diverges when t goes into the infinity. While let s=$e^{t}-1$, then by some simple calculation we know $\omega$(t) solves the normalized form equation is equivalent to $\tilde{\omega}$(s)=$e^{t}\omega$(t) solves the original K$\ddot{a}$hler-Ricci flow equation
                                       	
                                       	 $\frac{\partial}{\partial t}\omega$ = -Ric($\omega$).

                                       	 Therefore, the estimates we have proved before can be used in the proof of the original problem. Using the same method as the former sections, we now prove the existence of the K$\ddot{a}$hler- Einstein metric problem.
                                       	 
                                       	 $\mathbf{Theorem \; of \; existence \; of \; the \;K\ddot{a}hler-Einstein \; metric}$
                                       	 
                                       	 The solution to the equation

                                       	 $\frac{\partial\tilde{\mathit{g_{i\bar{j}}}}}{\partial t}$ = -$\tilde{\mathit{R_{i\bar{j}}}}$ - $\tilde{\mathit{g}}_{i\bar{j}}$,  $\tilde{\mathit{g_{i\bar{j}}}}$ = $\mathit{g_{i\bar{j}}}$ at t = 0 ,
                                       	 
                                       	 converges in C$^{\infty}$ to the unique K$\ddot{a}$hler-Einstein metric $\tilde{\mathit{g}}_{i\bar{j}}(\infty)$ which belongs to the negative first Chern class of M.
                                       	 
                                       	 $\mathit{Proof.}$
                                       	 
                                       	 First it's a K$\ddot{a}$hler Ricci flow equation which is a rescaling of 
                                       	 $\frac{\partial\tilde{\mathit{g_{i\bar{j}}}}}{\partial t}$ = -$\tilde{\mathit{R_{i\bar{j}}}}$ :
                                       	 if $u^{\prime}$(s) is a solution of the equation above, then u(t) = $\frac{u^{\prime}(s)}{s+1}$ for t = log(s+1), where s$\in$[0,$\infty$) , then  it has a short time solution, and by the Theorem of maximal existence time, it has an all time solution, this is due to the parabolic equation theory, and since we now have an zero order estimate of $\tilde{\mathit{g}}_{i\bar{j}}$ by Lemma 7.5 and Lemma 6.1, and an estimate for $\frac{\partial u}{\partial t}$ with its exponential decay by Lemma 6.1, then as the same arguement as the Proposition before, we use Schauder estimate and Interior regularity theory to obtain a C$^{\infty}$ estimate of $\tilde{\mathit{g}}_{i\bar{j}}$, then the C$^{0}$ estimate of u(t) gives the uniform C$^{\infty}$ estimate of it. While u(t) converges uniformly to u($\infty$) continuously when t goes into the infinity by Lemma 6.1,
                                       	 we prove that u(t) converges to u($\infty$) in the C$^{\infty}$ sense by contradiction. If there exists an integer k and $\epsilon$ positive, and a sequence $t_{i}$ goes into the infinity such that 
                                       	 
                                       	 $\mid\mid  u(t_{i})-u(\infty)\mid\mid_{C^{k}(M)}$ $\ge$ $\epsilon$ for any positive integer i,
                                       	 
                                       	 then because u($t_{i}$) has uniform C$^{k+1}$ bound then by Arzela-Ascoli Theorem, there exists a subsequence u($t_{i_{k}}$) converges to another limit, says $u^{\prime}(\infty)$ in the C$^{k}$ sense, but 
                                       	 
                                       	 $\mid\mid  u(t_{i})-u(\infty)\mid\mid_{C^{k}(M)}$ $\ge$ $\epsilon$ implies 
                                       	 $\mid\mid  u^{\prime}(\infty)-u(\infty)\mid\mid_{C^{k}(M)}$ $\ge$ $\epsilon$,

                                       	 so $u^{\prime}(\infty)\neq u(\infty)$ which is contradic to the uniqueness of the uniformly convergence of u(t). So we only get u(t) converges to u($\infty$) in the C$^{\infty}$ sense.
                                       	 Then this presents that $\frac{\partial u}{\partial t}$ converges to 0 as t goes into the infinity, because u converges to a constant u$_{\infty}$ therefore,from 
                                       	 $\tilde{\mathit{g_{i\bar{j}}}}$($\infty$) = $\mathit{g_{i\bar{j}}}$ + $\frac{\partial^{2}u(\infty)}{\partial\mathit{z^{i}}\partial\mathit{\bar{z}^{j}}}$,we get  $\tilde{\mathit{g}}_{i\bar{j}}$ converges to a constant of t when t goes into the infinity, 
                                       	 which means that $\frac{\partial \tilde{\mathit{g_{i\bar{j}}}}(t) }{\partial t}$ = 0,  then

                                       	 0 = -$\tilde{\mathit{R_{i\bar{j}}}}$($\infty$) - $\tilde{\mathit{g}}_{i\bar{j}}$($\infty$), 
                                       	 hence $\tilde{\mathit{R_{i\bar{j}}}}$($\infty$) =- $\tilde{\mathit{g}}_{i\bar{j}}$($\infty$), 
                                       	 
                                       	 so that's the K$\ddot{a}$hler-Einstein metric g we want.

                                       	 And the uniqueness follows. If $g^{\prime}$ is another K$\ddot{a}$hler-Einstein metric both belonging to the same negative first Chern class, then we can write: 
                                       	 
                                       	 $g^{\prime} = g+\frac{\sqrt{-1}}{2\pi}\partial\tilde{\partial}\varphi$, then by calculation Ric($g^{\prime}$) = Ric(g)-$\frac{\sqrt{-1}}{2\pi}\partial\tilde{\partial}\varphi$, so 
                                       	 
                                       	 $\frac{log det (g+\frac{\sqrt{-1}}{2\pi}\partial\tilde{\partial}\varphi)}{log det g}$ = $\varphi$ + C
                                       	 
                                       	 for C a constant,then by the maximum principle of the function $\varphi$ + C, the maximum and the minimum of $\varphi$ attains at the boundary , since by calculation as before arguement this $\varphi$+C satisfies the equation
                                       	 
                                       	 $\frac{\partial u}{\partial t}$ =   
                                       	 log det($\mathit{g^{\prime}}$ + $\frac{\partial^{2}u}{\partial\mathit{z^{i}}\partial\mathit{\bar{z}^{j}}}$ ) - log det($\mathit{g}$)  - u + f
                                       	  for $\varphi$+C= 0 when t=0,
                                       	 so we get the maximum and the minimum of the $\varphi$+C are  0, then the LHS is zero, then g$^{\prime}$=g. The proof of uniqueness actually comes from Calabi. Finally, the proof is completed, and we just justify the deformation method works well in K$\ddot{a}$hler-Einstein metric existence problem by some important estimates from the K$\ddot{a}$hler Ricci flow equation.

                                 $ \mathbf{References}$

                                 1.Cao, H D. Deformation of K$\ddot{a}$hler metrics to K$\ddot{a}$hler-Einstein metrics on compact Kahler manifolds. United States: N. p., 1986.
                                 
                                 2.
                                 Chern, Shiing-shen. “Characteristic Classes of Hermitian Manifolds.” Annals of mathematics 47.1 (1946): 85–121. Web.
                                 
                                 3.Richard S. Hamilton "Three-manifolds with positive Ricci curvature," Journal of Differential Geometry, J. Differential Geom. 17(2), 255-306, (1982)
                                 
                                 4.Yau, S.-T. (1978), On the ricci curvature of a compact K$\ddot{a}$hler manifold and the complex monge-ampére equation, I. Comm. Pure Appl. Math., 31: 339-411.
                                 
                                 5.
                                 Gilbarg, David, and Neil S. Trudinger. Elliptic Partial Differential Equations of Second Order. 1st ed. 1977. Berlin, Germany, Springer-Verlag, 1977. Web.
                                 
                                 6.
                                 Petersen, Peter. Riemannian Geometry. Third edition. Cham, Switzerland: Springer, 2016. Web.
                                 
                                 7.
                                 Carmo, Manfredo Perdigao do. Riemannian Geometry, 2008. Print.

                                 8.
                                 Elliptic Partial Differential Equations: Second Edition (Courant Lecture Notes) 2nd Edition. Han Qing, Lin Fang Hua,American Mathematical Society, 2011. Print.
                                 
                                 9.
                                 Jian Song, and Ben Weinkove. “Lecture Notes on the K$\ddot{a}$hler-Ricci Flow.” 	arXiv:1212.3653 [math.DG]
                                 
                                 10.
                                 PDE II Schauder estimate, Robert  Hasslhofer. Lecture notes, Toronto university https://www.math.toronto.edu/roberth/pde2/schauder$\_$estimates.pdf
                                 
                                 11.
                                 “The Ricci Flow; Techniques and Applications, Pt.1: Geometric Aspects.” SciTech Book News 31.2 (2007): n. pag. Print.
                                 
                                 12.Calabi, Eugenio. "On K$\ddot{a}$hler Manifolds with Vanishing Canonical Class". Algebraic Geometry and Topology: A Symposium in Honor of Solomon Lefschetz, edited by Ralph Hartzler Fox, Princeton: Princeton University Press, 2015, pp. 78-89. https://doi.org/10.1515/9781400879915-006
                                 
                                 13.Ben Weinkove, "The K$\ddot{a}$hler-Ricci flow on compact K$\ddot{a}$hler manifolds",	arXiv:1502.06855 [math.DG]
                                 
                                 14.Evans, Lawrence C.. “Partial Differential Equations, Second edition.” (2010).

\end{document}